\documentclass[11pt,a4paper]{article}

\usepackage{natbib}

\usepackage{algorithm}
\usepackage{algpseudocode}

\usepackage{amsmath,amsthm}
\usepackage{amssymb}

\usepackage[latin1]{inputenc}
\usepackage[T1]{fontenc}

\usepackage{multirow}

\usepackage[none]{hyphenat}

\usepackage{booktabs}

\usepackage{tikz}
\usetikzlibrary{calc}
\usetikzlibrary{positioning}
\usetikzlibrary{shapes.arrows}

\usepackage{caption}
\usepackage{pgfplots}

\newcommand{\boldm}[1] {\mathversion{bold}#1\mathversion{normal}}
\newcommand{\round}[1]{\ensuremath{\lfloor#1\rceil}}

\usepackage{color, colortbl}

\definecolor{Gray}{gray}{0.9}
\usepackage{setspace} 

\newcommand{\specialcell}[2][c]{%
	\begin{tabular}[#1]{@{}c@{}}#2\end{tabular}}

\usepackage[margin=2.3cm]{geometry}
\usepackage{xcolor}
\definecolor{myblue}{HTML}{D6EAF8} 

\tikzset{
	mybox/.style={rectangle,
		draw,
		fill= myblue,
		rounded corners,
		minimum width=2cm,
		inner sep=5pt,
		align=left,
		minimum height=1cm
	},
	myarrow/.style={draw=black,
		fill=white,
		minimum width=0.6cm,
		single arrow
	},
	downarrow/.style={draw=black,
	fill=white,
	minimum width=0.6cm,
	single arrow
	},
	longarrow/.style={draw=none,
		shading=axis,
		left color=white,
		right color=myblue,
		minimum width=0.6cm,
		single arrow,
		anchor=east
	}
}

\begin{document}


	\title{Air cargo load and route planning \\ in pickup and delivery operations}
	\author{A.C.P.~Mesquita (celio@ita.br) and C.A.A.~Sanches (alonso@ita.br) \\ Instituto Tecnol\'{o}gico de Aeron\'{a}utica - DCTA/ITA/IEC \\ Pra\c{c}a Mal. Eduardo Gomes, 50\\
		S\~{a}o Jos\'{e} dos Campos - SP - 12.228-900 - Brazil}

\maketitle

\begin{abstract}
In the aerial pickup and delivery of goods in a distribution network, transport aviation faces risks of load imbalance due to the urgency required for loading, immediate take-off, and mission accomplishment. Transport planners deal with trip itineraries, prioritisation of items, building up pallets, and balanced loading, but there are no commercially available systems that can integrally assist in all these requirements. This enables other risks, such as improper delivery, excessive fuel burn, and possible safety issues due to cargo imbalance, as well as a longer than necessary turn-around time. This {\it NP-hard}\/ problem, named {\it Air Cargo Load Planning with Routing, Pickup, and Delivery Problem (ACLP+RPDP)}, is mathematically modelled using standardised pallets in fixed positions. We developed a strategy to solve this problem, considering historical transport data from some Brazilian hub networks, and performed several experiments with a commercial solver, five known meta-heuristics, and a new heuristic designed specifically for this problem. By using a portable computer, our strategy quickly found practical solutions to a wide range of real problems in much less than operationally acceptable time.
\end{abstract}



\section{Introduction}\label{sec:Introduction}

The aviation industry adapts during global crises to keep supply chains moving. Air cargo provided complex expertise and the ability to access diverse destinations, delivering essential goods such as medicines, vaccine supplies, testing kits and other necessities with exceptional speed. This mode of transportation has become a preferred choice for governments, corporations and global companies in urgent need of transportation solutions.

Air cargo services are specially designed for organisations that require customised transportation, handle sensitive goods, or serve remote locations with limited routes. Air carriers typically use high-capacity cargo planes for economies of scale. Many cargo airlines have worldwide networks spread across destinations around the world.

A few years ago, \cite{BrandtStefan2019} defined the {\it Air Cargo Load Planning Problem} (ACLPP) as four sub-problems: {\it Aircraft Configuration Problem} (ACP), {\it Build-up Scheduling Problem} (BSP), {\it Air Cargo Palletization Problem} (APP), and {\it Weight and Balance Problem} (WBP). Several aspects were considered: item characteristics to be transported (dimensions, scores, dangerousness, etc.); types and quantities of {\it unit load devices} (ULDs); when these pallets are assembled; how items are allocated to pallets; in which positions these pallets are to be placed; how total cargo weight is balanced; etc.

However, it is crucial to highlight that there are still other important challenges in air cargo transport that go beyond the definition of ACLPP, especially with regard to routes, and pickup and delivery at each destination. In this context, at least two more important sub-problems can be considered: simultaneous pickup and delivery at each node, called the {\it Pickup and Delivery Problem} (PDP), and searching for the best benefit-cost route, which is a special case of the {\it Travelling Salesman Problem} (TSP).

Inefficient air transport plans can lead to unnecessary costs, extra routes, longer distances, and incorrect destinations. Unbalanced cargo increases fuel consumption due to altered aircraft pitch angles, increasing the risk of weight- and balance-related accidents. Balancing cargo is crucial for safe aerial transportation, as an improperly positioned centre of gravity (CG) can result in dangerous take-off and landing conditions and stall recovery issues. Despite technological advancements, many airlines still rely on manual aircraft loading and balancing, which can lead to flight delays. A rational decision-making process is essential to avoid creating inefficient or unsafe transport plans, considering high costs of fuel, maintenance, operation, outsourcing expenses, potential operational impairments, and safety risks due to unbalanced cargo. Solving this problem is vital for optimizing strategic scores, saving time and effort in loading, ensuring safety and balance, correct pickups and deliveries, and finding the best routes considering fuel use due to potential cargo imbalance.

The problem addressed in this work is the definition of the route of an aircraft that loads and unloads hundreds of items in several hubs, dealing with weight, volume, and balance constraints, and maximising the benefit-cost ratio. This challenge has become even more acute during the COVID-19 pandemic due to the need for urgent medical supplies. Delays in shipping essential items such as respirators to critical areas highlighted the need for optimised solutions.

This problem is extremely complex, as it deals with different objectives: defining a route that visits all hubs; maximising the items transported along this route, prioritising the essential ones; making sure items reach the correct destinations; ensuring aircraft safety constraints; and saving fuel so that the flight is sustainable. 

We have developed a heuristic process that can be run on a handheld computer, quickly providing a good solution to real instances of this problem. Solutions consist of flight itineraries, pickup and delivery plans, and the allocation of items onto pallets, ensuring the load-balancing constraints. Our method also reduces the stresses that transport planners are subject to, as they have to deal with extensive information in a short time frame.

To the best of our knowledge, this is the first time that an air cargo problem involving simultaneously APP, WBP, PDP, and TSP has been addressed. This new problem is named {\it Air Cargo Load Planning with Routing, Pickup, and Delivery Problem} (ACLP+RPDP). As we will describe in our mathematical modelling, these four sub-problems appear in an interconnected way in ACLP+RPDP and therefore cannot be solved independently.

As a real case study, we consider a crucial network for the {\it Brazilian Air Force}, as can be seen in Table \ref{tab:costs} and Figure \ref{fig:nodes}. Although there are other airports of interest, these nodes were chosen due to their high demand. Other Brazilian airports tend to have smaller transport requests, which are generally met less expensively by cabotage, rail, or road transport.

\begin{table}[htp]
	
	\begin{minipage}{0.58\linewidth}
		
		\caption{Distances between some Brazilian airports ($km$)} \label{tab:costs}
		\centering
		\footnotesize
		
		\newcolumntype{C}{>{\centering\arraybackslash}p{0.07\textwidth}}
		\newcolumntype{R}{>{\raggedleft\arraybackslash}p{0.07\textwidth}}			
		
		\begin{tabular}{C R R R R R R R}
			\toprule
			IATA*   & GRU   & GIG   & SSA   & CNF   & CWB   & BSB   & REC \\	
			\midrule	
			GRU     & 0	    &343	&1,439   &504    &358    &866    &2,114\\
			GIG	    & 343	&0	    &1,218   &371    &677    &935    &1,876\\
			SSA	    & 1,439	&1,218	&0	    &938    &1,788   &1,062   &676\\
			CNF	    & 504	&371	&938	&0	    &851    &606    &1,613\\
			CWB	    & 358	&677	&1,788	&851	&0	    &1,084   &2,462\\
			BSB	    & 866	&935	&1,062	&606	&1,084	&0	    &1,658\\
			REC	    & 2,114	&1,876	&676	&1,613	&2,462	&1,658	&0\\
			\bottomrule
			\multicolumn{8}{c}{*International Air Transport Association}\\
			\multicolumn{8}{c}{\small\textsuperscript{Source: www.airportdistancecalculator.com}}\\
		\end{tabular}
		\normalsize
		
	\end{minipage}\hfill 
	\begin{minipage}{0.42\linewidth}
		\centering
		\includegraphics[scale=0.21]{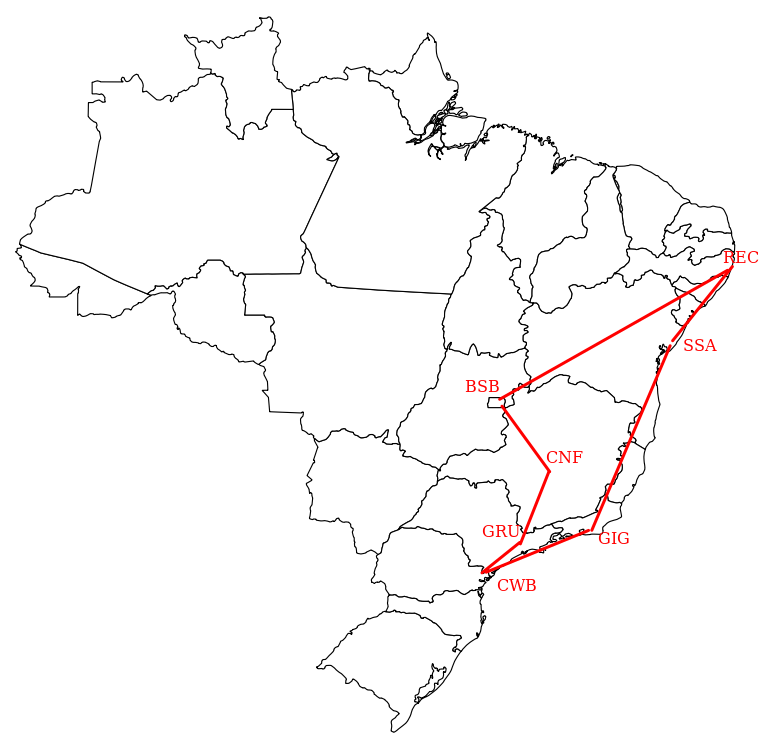}
		\captionof{figure}{A route with 7 airports}
		\label{fig:nodes}
	\end{minipage}
\end{table}

This article is organised into six more sections. In Section \ref{literature}, we make the literature review. In Section \ref{assumptions}, we present the context and requirements of ACLP+RPDP. In Section \ref{modelling}, we describe its mathematical modelling. In Section \ref{algorithms}, we describe the developed algorithms, whose results are presented in Section \ref{results}. Finally, our conclusions are in Section \ref{conclusions}.

\section{Literature review}
\label{literature}

The vast majority of operational research applied to air cargo is focused on challenges related to WBP, that is, the distribution of items on pallets to ensure load balancing. We can mention: \cite{LarsenMikkelsen1979}; \cite{Brosh1981}; \cite{Kevin1992}; \cite{Heidelberg1998}; \cite{fok2004optimizing}; \cite{KaluznyBohdanL2009Oalb}; \cite{Verstichel2011}; \cite{Limbourg2012}; \cite{RoesenerBarnes2016}; \cite{YangLiuGao2018}; \cite{zhao2021}; \cite{MiguelOptimalAP}.

Other authors have addressed pallet assembly (APP) on aircraft, possibly also considering load balancing (WBP): \cite{MongeauBes2003}; \cite{Chan2006}; \cite{RoesenerHall2014}; \cite{Vancroonemburg2014}; \cite{Paquay2016}; \cite{Paquay2018}; \cite{wong2020}; \cite{eugene2021}; \cite{zhao2023}.

In all these works, there is a great diversity of scenarios and solutions: some consider items in two dimensions, and others in three dimensions; some used integer programming, and others developed specific heuristics. Among the most recent, we can highlight:
\begin{itemize}
	
	\item \cite{RoesenerBarnes2016} proposed a heuristic to solve the {\it Dynamic Airlift Loading Problem} (DALP). Given a set of palletized cargo items that require transport between two nodes in a given time frame, the objective of this problem is to select an efficient subset of aircraft, partition the pallets into aircraft loads, and assign them to allowable positions in those aircraft.
	
	\item \cite{Paquay2016} presented a mathematical model to optimise the loading of heterogeneous 3D boxes on pallets with a truncated parallelepiped format. Its objective is to maximise the volume used in containers, considering load balancing constraints, the presence of fragile items, and the possibility of rotating these boxes. \cite{Paquay2018} developed some heuristics to solve this problem.
	
	\item \cite{YangLiuGao2018} modelled the air transport problem as a 2D packing problem and presented a heuristic for its optimisation in several aircraft, considering load balancing to minimise fuel consumption.
	
	\item \cite{wong2020} developed a mathematical model and a tool based on mixed integer programming for optimising cargo in aircraft with different pallet configurations. Balance constraints and the presence of dangerous items were considered. \cite{eugene2021} integrated this tool into a digital simulation model with a visualisation and validation system based on sensors that alert about load deviations.
	
	\item \cite{zhao2021} proposed a model for WBP based on Mixed Integer Programming (MIP). Instead of focusing on the CG deviation, the authors consider the original CG envelope of the aircraft, with a linearization method for its non-linear constraints.
	
	\item \cite{zhao2023} presented three models that use integer programming for air cargo planning and weight balance optimization: bi-objective optimization (BOM), combinatorial optimization (COM), and enhanced combinatorial optimization (IOM). Considering a Boeing 777F in several scenarios, the tests revealed performance problems: BOM is fast, but produces large CG deviation; COM offers accurate optimization, but with impractical runtimes; IOM provides a balanced solution, improving speed over COM, but requiring high computational demands in some cases. Although IOM stands out for its effectiveness, all models face trade-offs between speed, accuracy and computational efficiency. This work alerted us to potential performance issues in solution methods.
	
	\item \cite{MiguelOptimalAP} emphasized the proper selection of an aircraft, balancing factors such as cost, efficiency and limitations. These authors proposed a fuzzy linear programming model, which allows airlines to consider multiple objectives: maximizing payload, prioritizing specific items and minimizing operational costs. However, this approach does not take into account palletizing constraints.
	
\end{itemize}

On the other hand, there are several works that address PDP \citep{Golestanian, Meng, Bertsimas} or TSP \citep{Debnath, Cheikh, Ahmad, Xie} for unmanned aerial vehicles or aircraft, but none of them deal with APP and WBP.

\cite{LurkinSchyns2015} is the only work that simultaneously addresses an air cargo (WBP) and a flight itinerary (PDP) sub-problem. The authors demonstrated that this problem is {\it NP-hard}. Although it is innovative, strong simplifications were imposed by these authors: in relation to loading, APP was ignored; regarding routing, it is assumed that a predefined tour plan is restricted to only two legs. Referring directly to this work, \cite{BrandtStefan2019} comment: {\it However, not even these sub-problems are acceptably solved for real-world problem sizes, or models omit some practically relevant constraints}.

Table \ref{tab:sa} lists the literature on air cargo transport with the sub-problems involved. We also indicate whether the dimensions of the items were considered ({\bf 2D} or {\bf 3D}) and which solution method was used: heuristic search methods ({\bf Heu}), integer programming ({\bf Int}), or linear programming ({\bf Lin}).

\begin{table}[htp]
	\centering
	\caption{Air cargo transport: literature, sub-problems and features}  \label{tab:sa}
	\footnotesize
	\begin{tabular}{r|cccc|ccccc}
		\toprule
		Work & {\bf APP}  & {\bf WBP}  &  {\bf PDP}   &{\bf TSP}   & {\bf 2D}  & {\bf 3D}  & {\bf Heu}  & {\bf Int}  & {\bf Lin} \\
		\midrule
		\cite{LarsenMikkelsen1979}  & $.$        & $\bigstar$ & $.$          & $.$        & $.$       & $.$       & $\bigstar$ & $.$        &  $.$ \\
		\cite{Brosh1981}  & $.$ & $\bigstar$  & $.$   & $.$ & $.$ & $.$   & $.$  & $.$  &  $\bigstar$ \\
		\cite{Kevin1992}  & $.$ & $\bigstar$  & $.$ & $.$ & $.$ &$.$   & $.$  & $\bigstar$  &  $.$ \\
		\cite{Heidelberg1998}  & $.$ & $\bigstar$  & $.$ & $.$ & $\bigstar$ &$.$   & $\bigstar$  & $.$  &  $.$ \\
		\cite{MongeauBes2003}    & $\bigstar$ & $\bigstar$   & $.$ & $.$ & $.$ & $.$   & $.$  & $\bigstar$  &  $.$ \\
		\cite{fok2004optimizing} & $.$ & $\bigstar$   & $.$ & $.$ & $.$ & $.$   & $.$  & $\bigstar$  &  $.$ \\	
		\cite{Chan2006}  & $\bigstar$ & $.$    & $.$ & $.$ & $.$ & $\bigstar$  & $\bigstar$  & $.$  &  $.$ \\
		\cite{KaluznyBohdanL2009Oalb}  & $.$ & $\bigstar$  & $.$  & $.$ & $\bigstar$ &$.$  & $.$  & $\bigstar$  &  $.$ \\
		\cite{Verstichel2011}   & $.$ & $\bigstar$    & $.$ & $.$ & $.$ & $.$   & $.$  & $\bigstar$  &  $.$ \\	
		\cite{Limbourg2012} & $.$ & $\bigstar$  & $.$ & $.$ & $.$ & $.$   & $.$  & $\bigstar$  &  $.$ \\
		\cite{RoesenerHall2014}  & $\bigstar$ & $\bigstar$  & $.$  & $.$ & $.$ & $\bigstar$   & $.$  & $\bigstar$  &  $.$ \\
		\cite{Vancroonemburg2014}  & $\bigstar$ & $\bigstar$   & $.$ & $.$ & $.$ & $.$   & $.$  & $\bigstar$  &  $.$ \\
		\cite{LurkinSchyns2015} & $.$ & $\bigstar$  & $\bigstar$ & $.$  & $.$ & $.$   & $.$  & $\bigstar$  &  $.$ \\
		\cite{RoesenerBarnes2016}  & $.$ & $\bigstar$   & $.$ & $.$ & $.$ & $.$   & $\bigstar$  & $.$  &  $.$ \\
		\cite{Paquay2016,Paquay2018}  & $\bigstar$ & $\bigstar$ & $.$ & $.$ & $.$ & $\bigstar$ & $\bigstar$  & $\bigstar$ & $.$ \\
		\cite{YangLiuGao2018} & $.$ & $\bigstar$  & $.$  & $.$ & $\bigstar$  & $.$ & $\bigstar$ & $.$  & $.$ \\
		\cite{wong2020} & $\bigstar$  & $\bigstar$  & $.$  & $.$   & $.$  & $.$ & $.$ & $\bigstar$  & $.$  \\
		\cite{eugene2021} & $\bigstar$ & $\bigstar$ & $.$  & $.$   & $.$ & $.$ & $.$ & $\bigstar$  & $.$  \\
		\cite{zhao2021} & $.$ & $\bigstar$ & $.$  & $.$  & $.$ & $.$ & $.$  & $\bigstar$ &  $.$ \\
		\cite{zhao2023} & $\bigstar$ & $\bigstar$ & $.$  & $.$  & $.$ & $.$ & $.$  & $\bigstar$ &  $.$  \\
		\cite{MiguelOptimalAP} & $.$ & $\bigstar$ & $.$  & $.$  & $.$ & $.$ & $.$  & $\bigstar$ &  $.$  \\
		{\bf This article}   & $\bigstar$ & $\bigstar$  & $\bigstar$& $\bigstar$ & $.$ & $.$ & $\bigstar$ & $\bigstar$   &  $.$  \\
		\bottomrule 
	\end{tabular}
	\normalsize 
\end{table}

As can be seen, none of these papers address air cargo palletization and load balancing with route optimisation in a multi-leg transport plan for a single aircraft. Our work is the first to address a real air transport problem in which APP, WBP, PDP and TSP arise in an interconnected way.

\section{Context and assumptions}
\label{assumptions}

In this section, we describe the context of the problem addressed in this work as well as the assumptions considered.

\subsection{Operational premises}

As we are dealing with an extremely complex and diverse problem, we decided to establish some simplifying characteristics:

\begin{itemize}
	\setlength\itemsep{-0.3em}
	
	\item At each node of the tour, the items to be allocated are characterized by weight, volume, scores, and previously known destinations. We leave the consideration of 2D or 3D items to a future work.
	
	\item We considered a unique pallet type: the {\it 463L Master Pallet}, a common size platform for bundling and moving air cargo. It is the primary air cargo pallet for more than 70 Air Forces and many air transport companies. This pallet has a capacity of $4,500 kg$ and $13.7 m^3$, which may be limited by its position along the cargo bay. It is equipped for locking into cargo aircraft rail systems, and includes tie-down rings to secure nets and cargo loads, which in total weighs $140 kg$. For more information, see {\tt www.463LPallet.com}.
	
	\item All items allocated on a pallet must have the same destination. A pallet which has not yet reached its destination may receive more items, although it is known that these operations of removing restraining nets increase handling time and the risk of improper delivery. We do not consider oversized cargo in this work, but only cargo items that fit on these pallets.
	
	\item Finally, as we are interested in minimizing fuel costs, we disregarded others costs not directly associated with aircraft flight, such as handling.
	
\end{itemize}

\begin{figure}[H]
	\centering
	\includegraphics[scale=0.3]{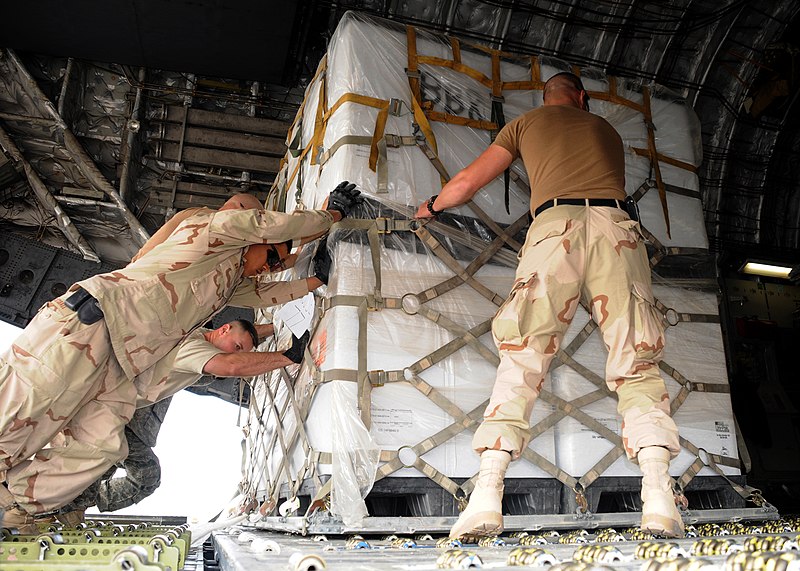}
	\caption{A {\it packed content}\/ on 463L pallet inside a {\it Boeing C-17}}
	\small\textsuperscript{Source: From Wikimedia Commons, the free media repository}
	\label{fig:packed}
\end{figure}

Throughout this text, we call {\it packed content}\/ (see Figure \ref{fig:packed}) a set of items of the same destination stacked on a pallet and covered with a restraining net. It is considered a single item, having the same attributes as its components, whose values are the sum of individual scores, weights, and volumes. To ensure accuracy in pickup and delivery operations, packed content must remain on board until its destination.

\subsection{Aircraft parameters and load balancing}

We consider real-world scenarios, where Table \ref{tab:larger} shows the aircraft parameters. $p_i$\/ are pallets, $1\leq i \leq 18$, whose weight and volume limits are $W_i$\/ and $V_i$, respectively. $D_i^{long}$\/ and $D_i^{lat}$\/ are, respectively, the longitudinal and lateral distances of each pallet centroids to the aircraft CG along both axes. These distances will be used in the calculation of the torque, referring to the items allocated on each pallet. In this aircraft, as the ramp has an inclination of $25^{\circ}$, we made the necessary corrections in $D_i^{long}$, $W_i$\/ and $V_i$\/ of the corresponding pallets ($p_1$, $p_2$, $p_3$, and $p_4$).

\begin{table}[htp]
	\centering
	\caption{Aircraft parameters}  \label{tab:larger}
	\footnotesize
	\begin{tabular}{c | c c c c c c c c c}
		\toprule
		& \multicolumn{3}{c}{$Payload$: 75,000kg} & \multicolumn{3}{c}{$limit^{CG}_{long}$: $1.170m$} &
		\multicolumn{3}{c}{$limit^{CG}_{lat}$: $0.19m$} \\
		\midrule
		\multirow{2}{*}{{\boldm{$p_i$}}}  & $p_{17}$ & $p_{15}$ & $p_{13}$ & $p_{11}$ & $p_{9}$ & $p_{7}$ & $p_{5}$ & $p_{3}$ & $p_{1}$ \\
		& $p_{18}$ & $p_{16}$ & $p_{14}$ & $p_{12}$ & $p_{10}$ & $p_{8}$ & $p_{6}$ & $p_{4}$ & $p_{2}$ \\
		\midrule 
		\multirow{2}{*}{\boldm{$D_i^{long}$} ($m$)} & -17.57 & -13.17 & -8.77 & -4.40 & 0 & 4.40 & 8.77 & 11.47 & 14.89 \\
		& -17.57 & -13.17 & -8.77 & -4.40 & 0 & 4.40 & 8.77 & 11.47 & 14.89 \\			
		\midrule 
		\multirow{2}{*}{\boldm{$D_i^{lat}$} ($m$)}  & 1.32 & 1.32 & 1.32 & 1.32 & 1.32 & 1.32 & 1.32 & 1.32 & 1.32 \\
		& -1.32 & -1.32 & -1.32 & -1.32 & -1.32 & -1.32 & -1.32 & -1.32 & -1.32 \\	
		\midrule
		{\boldm{$W_i$}} ($kg$)      &   4,500   &    4,500  &   4,500   &  4,500    & 4,500     & 4,500     & 4,500     & 3,000    & 3,000   \\
		{\boldm{$V_i$}} ($m^3$)   &   14.8   &   14.8   &  14.8    &  14.8    & 14.8     & 14.8     & 14.8     & 10.0    & 7.0 \\	
		\midrule	
		\textbf{Fuel cost}  & \multicolumn{9}{c}{ $c_d$ = US\$ $4.90/km$ }	 \\
		\midrule
		\textbf{Fuel consumption rate}  & \multicolumn{9}{c}{ $c_g = 5\%$} \\		
		\midrule
		\textbf{Maximum weight}  & \multicolumn{9}{c}{ $W_{max} = \sum_i W_i = 75,000kg$} \\
		\bottomrule
	\end{tabular}
	\normalsize 
\end{table}

This aircraft spends $c_d$\/ dollars per kilometre flown and can carry up to $W_{max}$\/ of cargo distributed on the pallets. The fuel penalty $c_g$\/ is the percentage of cost increase due to the CG deviation on the longitudinal axis, estimated at 5.0\%. It is important to consider that $c_g$\/ tends to zero as the aircraft attitude tends to be level. As the CG deviation varies from $0$\/ to $limit^{CG}_{long}$, this fuel penalty varies from $0$\/ to $c_g$.

\begin{figure}[htp]
	\centering
	\includegraphics[scale=0.22]{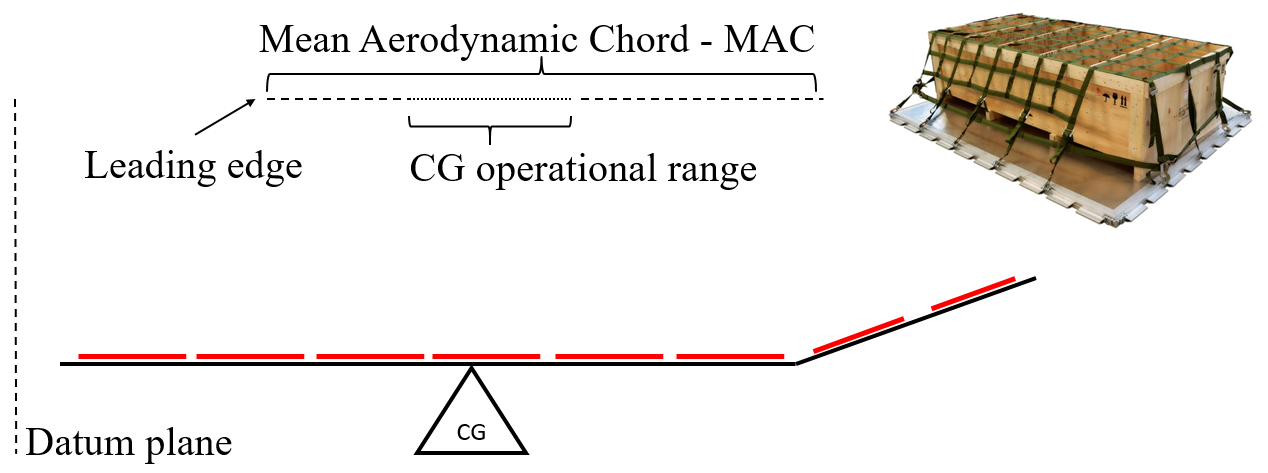}
	\caption{Aircraft longitudinal cut, where red lines are pallets positions}
	\label{fig:lateral}
\end{figure}

The torque applied to the aircraft must keep its CG in the operational range, which corresponds to a fixed percentage of the {\it Mean Aerodynamic Chord} \footnote{Chord is the distance between the leading and trailing edges of the wing, measured parallel to the normal airflow over the wing. The average length of the chord is known as the {\it Mean Aerodynamic Chord} (MAC).} which is considered $1.17m$ for the aircraft of this work (see Figure \ref{fig:lateral}).

We also make the following assumptions:

\begin{itemize}
	\setlength\itemsep{-0.3em}
	\item on each pallet, the items are distributed in such a way that their CG coincides with the centroid of the pallet, because builders are well-trained to do so;
	\item the CG of the total load must be at a maximum longitudinal distance of $limit^{CG}_{long}$\/ from the CG of the aircraft;
	\item the CG of the total load must be at a maximum lateral distance of $limit^{CG}_{lat}$\/ from the CG of the aircraft;
	\item the pallets are distributed in two identical rows (with odd and even indices, respectively), and the centroid of $p_i$\/ is at a distance $D^{lat}_{i}$\/ from the centreline of the aircraft;
	\item when there are items or packed contents in $p_i$, the common destination of this load will be assigned to variable $T_i$. 
\end{itemize}

\section{The mathematical modelling}
\label{modelling}

In this section, we present the mathematical modelling of ACLP+RPDP in Tables \ref{tab:struc}, \ref{tab:graph}, \ref{tab:functions}, and \ref{tab:constraints}, with their corresponding descriptions.

In Table \ref{tab:struc}, we describe the problem structure: nodes and their permutations, distances and associated costs, pallets characteristics, items available for shipment at each node, and packed contents shipped. The item $j$\/ in node $k$\/ has score $s_j$, weight $w_j$, volume $v_j$, and destination $to_j \in L_k$. Similarly, the packed content $q$, that remains on board at node $k$, has total weight $w_q$, total volume $v_q$, and destination $to_q \in L_k$. Packed contents that were destined to node $k$\/ are unloaded when the aircraft arrives there; that is, they are not considered in $Q_k$.

\begin{table}[htp]
	\centering
	\footnotesize
	\caption{Problem structure}  \label{tab:struc}
	\begin{tabular}{l l}
		\toprule
		\textbf{Notation}	& \textbf{Description} \\
		\midrule
		$L = \{ 0, 1, \ldots, K \}$ & Set of $K+1$\/ nodes of the tour, where node 0 is the base \\
		\midrule
		$\pi$               & A permutation between nodes 1, \ldots, $K$\\
		\midrule
		$S_K$               & Set of $K!$\/ permutations \\
		\midrule
		$\pi_k$             & The $k^{th}$\/ node of tour $\pi$, $1 \leq k \leq K$\\
		\midrule
		Tour $\pi$          & $\{0,\pi_1,\ldots, \pi_K, 0\}$ \\
		& For ease of notation, $\pi_0 = \pi_{K+1} = 0$\\
		\midrule
		$L_k$               & Set of remaining nodes of tour at node $k$, $0 \leq k \leq K$\\
		& By definition, $L_0 = L$\\
		\midrule
		$d(a,b)$            & Distance from node $a$\/ to node $b$, where $0 \leq a, b \leq K$ \\
		&   By definition, $d(a,a) = 0, \forall a$\\
		\midrule
		$C=\left[c_{a,b}\right]$    & Cost matrix of flights, where $c_{a,b} = c_d \times d(a,b)$\\
		\midrule
		$M = \{1, \ldots, m \}$  & Set of $m$\/ pallets in specific positions within the aircraft \\
		&  See Table \ref{tab:larger}, where $m=18$ \\
		\midrule
		$N_k = \{1, \ldots, n_k \}$ & Set of $n_k$\/ items available for loading at node $k$, $1 \leq j \leq n_k$, $0 \leq k \leq K$\\
		\midrule
		$N = \bigcup_{0 \leq k \leq K} N_k$ & Set of items in all nodes along a tour\\
		\midrule
		$Q_k = \{1, \ldots, m_k \}$ & Set of $m_k \leq m$\/ packed contents at node $k$, $1 \leq q \leq m_k$, $0 \leq k \leq K$ \\ 
		& By definition, $m_0 = 0$ and $Q_0 = \varnothing$\\
		\bottomrule
	\end{tabular}
	\normalsize
\end{table}

Table \ref{tab:graph} contains decision variables and the ACLP+RPDP allocation graph.

\begin{table}[htp]
	\centering
	\footnotesize
	\caption{Decision variables and allocation graph}  \label{tab:graph}
	\begin{tabular}{l l}
		\toprule
		\textbf{Notation}	& \textbf{Description} \\
		\midrule
		$X_{ij}^{\pi_k}$\/ and $Y_{iq}^{\pi_k}$ & Binary variables, where $1 \leq i \leq m$, $1 \leq j \leq n_{\pi_k}$, $1 \leq q \leq m_{\pi_k}$\/ and $0 \leq k \leq K$\\
		\midrule
		$X_{ij}^{\pi_k} = 1$ & If item $j$\/ at node ${\pi_k}$\/ is assigned to pallet $i$, and 0 otherwise\\
		\midrule
		$Y_{iq}^{\pi_k} = 1$ & If packed content $q$\/ at node ${\pi_k}$\/ is assigned to pallet $i$, and 0 otherwise\\
		\midrule
		$T_i^{\pi_k} \in L_{\pi_k}$ &  Destination of items and packed contents assigned to pallet $i$\ at node ${\pi_k}$\\
		\midrule
		$G_{\pi_k}(V_{\pi_k}, E_{\pi_k})$ & Allocation graph at node ${\pi_k}$\\
		\midrule
		$V_{\pi_k} = M \cup N_{\pi_k} \cup Q_{\pi_k}$ & Allocation graph vertices at node ${\pi_k}$: pallets, items and packet contents\\
		\midrule
		$E_{N_{\pi_k}}$ & Allocation graph edges at node ${\pi_k}$, corresponding to shipped items\\
		\midrule
		$E_{Q_{\pi_k}}$ & Allocation graph edges at node ${\pi_k}$, corresponding to packed contents\\
		\midrule
		$E_{\pi_k} = E_{N_{\pi_k}} \cup E_{Q_{\pi_k}}$ & Allocation graph edges at node ${\pi_k}$\\
		\midrule
		$(i, j) \in E_{N_{\pi_k}}$ & If $X_{ij}^{\pi_k} = 1$, where $i$\/ is a pallet and $j$\/ is a item at node ${\pi_k}$ \\
		\midrule
		$(i, q) \in E_{Q_{\pi_k}}$ & If $Y_{iq}^{\pi_k} = 1$, where $i$\/ is a pallet and $q$\/ is a packed content at node ${\pi_k}$\\
		\bottomrule
	\end{tabular}
	\normalsize
\end{table}

The calculus functions of ACLP+RPDP are described in Table \ref{tab:functions}.

\begin{table}[htp]
	\centering
	\footnotesize
	\caption{Calculus functions} \label{tab:functions}
	\begin{tabular}{ c l }
		\toprule
		\textbf{Function} & \textbf{Description} \\
		\midrule
		(\ref{eq:scores}) & Total score of transported items throughout tour $\pi$\\
		\midrule
		(\ref{eq:tau})    & Longitudinal torque applied by loaded pallets at node ${\pi_k}$\\
		\midrule
		(\ref{eq:costs})  & Total cost of fuel on tour $\pi$\/ (distances and CG longitudinal deviations) \\
		\midrule
		(\ref{eq:departing}) & Set of not visited nodes at node ${\pi_k}$\\
		\midrule
		(\ref{eq:LatItem}) & Lateral torque at node ${\pi_k}$ (shipped items) \\
		\midrule
		(\ref{eq:LatPacked}) & Lateral torque at node ${\pi_k}$ (packed contents)\\
		\midrule
		(\ref{eq1}) & Objective function of ACLP+RPDP\\
		\bottomrule
	\end{tabular}
	\normalsize
\end{table}

\begin{equation} \label{eq:scores}
	\tilde{s}_\pi = \sum_{k=0}^{K} \sum_{i=1}^{m} \sum_{j=1}^{n_{\pi_k}} X_{ij}^{\pi_k} \times s_j
\end{equation}

\begin{equation} \label{eq:tau}
	\tau_{\pi_k} = \sum_{i=1}^{m} \Big [ D_i^{long} \times \Big ( \sum_{j=1}^{n_{\pi_k}} X_{ij}^{\pi_k} \times w_j +  \sum_{q=1}^{m_{\pi_k}} Y_{iq}^{\pi_k} \times w_q \Big ) \Big ] \Big / W_{max} \times limit^{CG}_{long}; \ k \in \{0, \ldots, K\}
\end{equation}

\begin{equation} \label{eq:costs}
	\tilde{c}_{\pi} = \sum_{k=0}^{K} \Big [ c_{\pi_k, \pi_{k+1}}\times(1+c_g\times|\tau_{\pi_k}|) \Big ] 
\end{equation}

\begin{equation} \label{eq:departing}
	L_{\pi_k} = L_{\pi_{k-1}} - \{\pi_k\}; \ k \in \{1, \ldots, K\}
\end{equation}

\begin{equation} \label{eq:LatItem}
	\epsilon_{\pi_k}^t = \sum_{i=1}^{m} \Big [ D_i^{lat} \times \sum_{j=1}^{n_{\pi_k}} \Big ( X_{ij}^{\pi_k} \times w_j \times (i\%2) - X_{ij}^{\pi_k} \times w_j \times (i+1)\%2 \Big ) \Big ] \Big / W_{max} \times limit^{CG}_{lat}
\end{equation}

\begin{equation} \label{eq:LatPacked}
	\epsilon_{\pi_k}^a = \sum_{i=1}^{m} \Big [ D_i^{lat} \times \sum_{q=1}^{m_{\pi_k}} \Big ( Y_{iq}^{\pi_k} \times w_q \times (i\%2) - Y_{iq}^{\pi_k} \times w_q \times (i+1)\%2 \Big ) \Big ] \Big / W_{max} \times limit^{CG}_{lat}
\end{equation}

\begin{equation} \label{eq1}
	\max_{\pi \in S_K} f_\pi = \tilde{s}_\pi/\tilde{c}_\pi
\end{equation}

Longitudinal (\ref{eq:tau}) and lateral torques (\ref{eq:LatItem}, \ref{eq:LatPacked}) are calculated in proportion to the highest torque supported by the aircraft. As there are two rows of pallets, one on each side of the centerline, we use the operator {\it modulo}\/ ($\%$) to calculate lateral torques. In our experiments, we found that the magnitude of these lateral torques was always minimal, so we decided to ignore them in the fuel consumption (\ref{eq:costs}). The objective of ACLP+RPDP (\ref{eq1}) is to find a permutation $\pi \in S_K$ with the corresponding allocation of items on pallets at each node that maximises the function $f_\pi = \tilde{s}_\pi/\tilde{c}_\pi$.

Finally, ACLP+RPDP constraints related to each node ${\pi_k}$\/ are described in Table \ref{tab:constraints}.

\begin{table}[htp]
	\centering
	\footnotesize
	\caption{Constraints} \label{tab:constraints}
	\begin{tabular}{ c l }
		\toprule
		\textbf{Constraint} & \textbf{Description} \\
		\midrule
		(\ref{eq:torqlong}, \ref{eq:torqlat}) & Longitudinal and lateral torques must be within aircraft limits \\
		\midrule
		(\ref{eq:app2}, \ref{eq:app3}) & Items allocated to each pallet cannot exceed its weight and volume limits\\
		\midrule
		(\ref{eq:app4}) & At most, each item is associated with a single pallet \\
		\midrule
		(\ref{eq:app5}) & Packed contents that have not yet reached their destination must remain on board\\
		\midrule
		(\ref{eq22}, \ref{eq23}) & Items allocated on the same pallet must have the same destinations\\
		\midrule
		(\ref{eq24}, \ref{eq25}) & If there is a packed content on the pallet, it must also have the same destination as other items\\
		
		\bottomrule
	\end{tabular}
	\normalsize
\end{table}

\begin{equation} \label{eq:torqlong}
	| \tau_{\pi_k} | \leq 1;\ k \in \{0, \ldots, K\}
\end{equation}

\begin{equation} \label{eq:torqlat}
	| \epsilon_{\pi_k}^t + \epsilon_{\pi_k}^a| \leq  1
\end{equation}

\begin{equation} \label{eq:app2}
	\sum_{j=1}^{n_{\pi_k}} X_{ij}^{\pi_k} \times w_j + \sum_{q=1}^{m_{\pi_k}} Y_{iq}^{\pi_k} \times w_q  \leq W_i; \ i \in \{1, \ldots, m\}
\end{equation}

\begin{equation} \label{eq:app3}
	\sum_{j=1}^{n_{\pi_k}} X_{ij}^{\pi_k} \times v_j + \sum_{q=1}^{m_{\pi_k}} Y_{iq}^{\pi_k} \times v_q  \leq\ V_i; \ i \in \{1, \ldots, m\}
\end{equation}

\begin{equation} \label{eq:app4}
	\sum_{i=1}^{m} X_{ij}^{\pi_k} \leq 1; \ j \in \{1, \ldots, n_{\pi_k}\}
\end{equation}

\begin{equation} \label{eq:app5}
	\sum_{i=1}^{m} Y_{iq}^{\pi_k} = 1;\ to_q \in L_{\pi_k}; \ q \in \{1, \ldots, m_{\pi_k}\}
\end{equation}

\begin{equation} \label{eq22}
	X_{ij}^{\pi_k} <= X_{ij}^{\pi_k} \times (T_i^{\pi_k} - to_j + 1); \ i \in \{1, \ldots, m\}; \ j \in \{1, \ldots, n_{\pi_k}\}
\end{equation}

\begin{equation} \label{eq23}
	X_{ij}^{\pi_k} <= X_{ij}^{\pi_k} \times (to_j - T_i^{\pi_k} + 1); \ i \in \{1, \ldots, m\}; \ j \in \{1, \ldots, n_{\pi_k}\}
\end{equation}

\begin{equation} \label{eq24}
	Y_{iq}^{\pi_k} <= Y_{iq}^{\pi_k} \times (T_i^{\pi_k} - to_q + 1); \ i \in \{1, \ldots, m\}; \ q \in \{1, \ldots, m_{\pi_k}\}
\end{equation}

\begin{equation} \label{eq25}
	Y_{iq}^{\pi_k} <= Y_{iq}^{\pi_k} \times (to_q - T_i^{\pi_k} + 1); \ i \in \{1, \ldots, m\}; \ q \in \{1, \ldots, m_{\pi_k}\}
\end{equation}

Once the assumptions and the mathematical modelling are presented, it is possible to see that ACLP+RPDP is {\it NP-hard}. In a similar way to \cite{LurkinSchyns2015}, consider the simple case where $K=1$\/ (one leg), $m=2$\/ (two pallets around the aircraft CG), $2n$\/ sufficiently light items with same scores in node 0, and no items in node 1. Under these conditions, through polynomial reductions for the {\it Set-Partition Problem}, it is possible to demonstrate that the decision problem associated with ACLP+RPDP is {\it NP-complete}. For more details, see \cite[p.~6]{LurkinSchyns2015}.

\section{Solution process}
\label{algorithms}

Throughout our research, we have thoughtfully described ACLP+RPDP in standard MIP format and found that no solver can handle its practical cases in a feasible time. Thus, as ACLP+RPDP is highly complex, involving four intractable and interconnected sub-problems, we decided to focus on {\it real cases}, developing quick node-by-node solutions, not necessarily optimal, but which would allow us to obtain a complete tour.

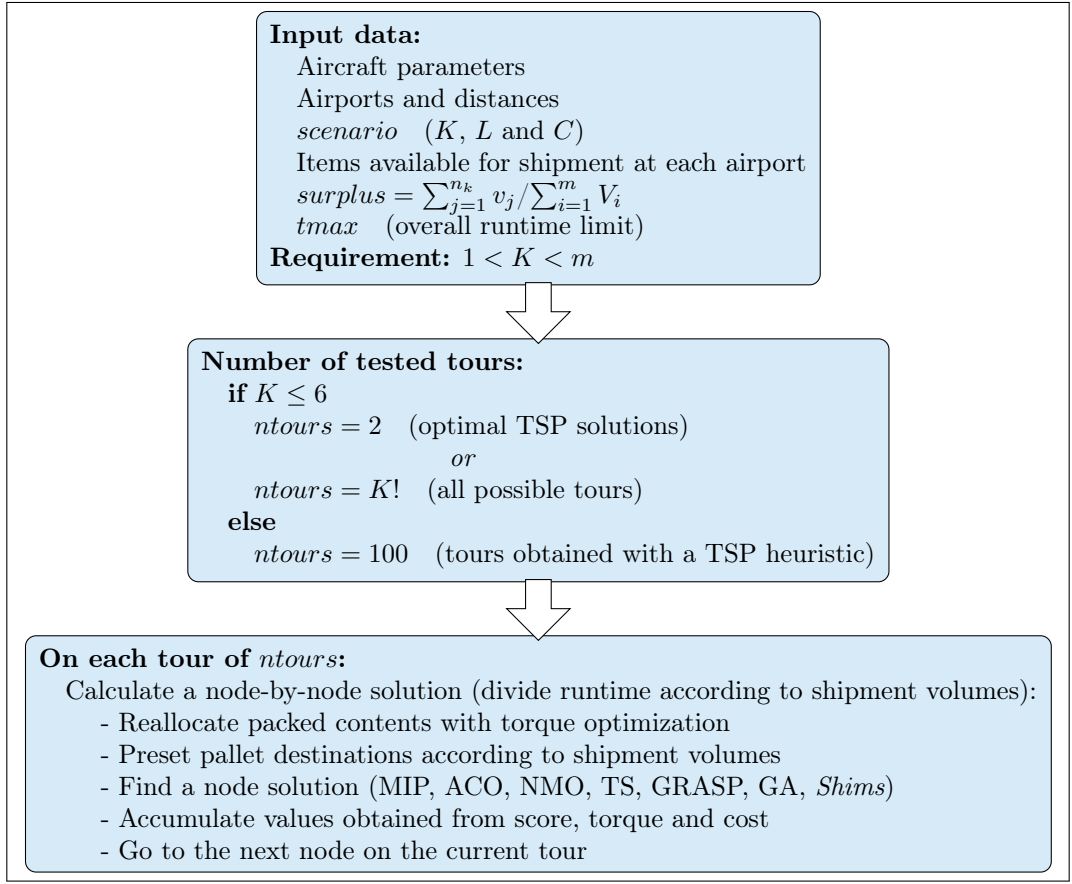
\begin{figure}[htp]
	\centering
	\fbox{ \small
		\begin{tikzpicture}
			\node [mybox] (one) { \small {\bf Input data:}\\
				\small \ \ \ {Aircraft parameters}\\
				\small \ \ \ {Airports and distances}\\
				\small \ \ \ {$scenario$ \ \ ($K$, $L$\/ and $C$)}\\
				\small \ \ \ {Items available for shipment at each airport}\\						
				\small \ \ \ {$surplus = \sum_{j=1}^{n_k} v_j$/$\sum_{i=1}^{m} V_i$ }\\		
				\small \ \ \ {$tmax$ \ \ (overall runtime limit)}\\
				\small {{\bf Requirement:} $1 < K < m$}
			};
			
			\node [mybox, below=0.7cm of one] (two) { \small {\bf Number of tested tours:}\\
				\small \ \ \ {{\bf if} $K \leq 6$}\\
				\small \ \ \ \ \ \ {$ntours = 2$ \ \ (optimal TSP solutions)}\\
				\small \ \ \ \ \ \ \ \ \ \ \ \ \ \ \ \ \ \ \ \ \ \ \ \ \ \ \ \ {\it or}\\				
				\small \ \ \ \ \ \ {$ntours = K!$ \ \ (all possible tours)}\\
				\small \ \ \ {\bf else}\\
				\small \ \ \ \ \ \ {$ntours = 100$ \ \ (tours obtained with a TSP heuristic)}
			};
			
			\node [mybox, below=0.7cm of two] (three) { \small {\bf On each tour of $ntours$:}\\
				\small \ \ \ {Calculate a node-by-node solution (divide runtime according to shipment volumes):}\\
				\small \ \ \ \ \ \ \ {- Reallocate packed contents with torque optimization}\\
				\small \ \ \ \ \ \ \ {- Preset pallet destinations according to shipment volumes}\\
				\small \ \ \ \ \ \ \ {- Find a node solution (MIP, ACO, NMO, TS, GRASP, GA, {\it Shims}) }\\
				\small \ \ \ \ \ \ \ {- Accumulate values obtained from score, torque and cost}\\
				\small \ \ \ \ \ \ \ {- Go to the next node on the current tour}
			};

			\node[myarrow, rotate=-90] at ([yshift=-8pt]two.south) {\phantom{\tiny{aaa}}};
			\node[myarrow, rotate=-90] at ([yshift=-8pt]one.south) {\phantom{\tiny{aaa}}};		
		\end{tikzpicture} 
	}
	\caption{Solution process}
	\label{fig:process}
\end{figure}

In practical cases, we know that a common aircraft has $m=18$\/ pallets, flight itineraries have $K \leq 6$\/ nodes plus the base, and each node has hundreds of items to be shipped. We also know that missions with fewer nodes are more frequent than longer ones. Under these circumstances, we can adopt some important strategies summarized in Figure \ref{fig:process}:

\begin{itemize}
	
	\item We consider that the number of destinations is smaller than the number of pallets ($K < m$), and we avoid the trivial case where $K=1$. With this premise, we can preset the destinations of the pallets at each shipping node, reserving a number of pallets proportional to the volume available for each destination. We could have used another criterion, but it was observed in the experiments that the volume is more constrictive in airlift.
	
	\item An important parameter is the number $ntours$\/ of tours tested. In practical cases where $K \leq 6$, we have the possibility to check all possible tours ($ntours = K!$). In this situation, as $K$\/ is small, we can also specially analyse the two optimal solutions of the corresponding TSP ($ntours = 2$). Finally, in cases where $K>6$, we will use a heuristic to select $100$\/ tours of low length ($ntours=100$), and search among them for the one that provides the best value for the objective function.
	
	\item To compare the performance of each strategy, an overall runtime limit $tmax$\/ is established and divided by $ntours$\/ tours. In turn, the runtime limit for each tour will be distributed among its nodes in proportion to the volume available for boarding.
	
	\item At each node of a tour, the packed contents that remain on board are reallocated on pallets in order to minimize torque on the aircraft. This calculation is done quickly using a MIP solver. Then, the destinations of the pallets are previously defined in proportion to the shipment volume. Finally, considering the runtime limit of each node, we will use a MIP solver and five well-known meta-heuristics to find the best allocation of shipping items: {\it Ant Colony Optimization} (ACO), {\it Noising Method Optimization} (NMO), {\it Tabu Search} (TS), {\it Greedy Randomized Adaptive Search Procedure} (GRASP), and {\it Genetic Algorithm} (GA). We will also introduce a very fast heuristic, developed specifically for this problem, called {\it Shims}.
	
	\item We will generate benchmarks using the $surplus$\/ parameter, which is a value in $\{1.2, 1.5, 2.0\}$. It corresponds, at each node $k$, to the ratio between the sum of the volumes of the items and the load capacity of the pallets ($surplus = \sum_{j=1}^{n_k} v_j$/$\sum_{i=1}^{m} V_i$). This parameter allows us to verify the different behaviour of each method, according to $scenario$\/ and the quantity of items available for shipment.
	
	\item We will do tests by varying the number $K$\/ of destinations, the set $L$\/ of nodes, and the costs $C$. Each group of values tested is called $scenario$, according to Tables \ref{tab:costs} and \ref{tab:scenarios}, where $1<K \leq 6$. After finding the method with the best performance in node-by-node solution, we will test it in solving cases with $K > 6$.
	
\end{itemize}

Algorithm \ref{alg:main} is the main program of this solution process, which parameters are $scenario$, $surplus$, $tmax$, and $ntours$. The input data is obtained from Tables \ref{tab:costs}, \ref{tab:larger} and \ref{tab:scenarios} (lines \ref{main:M}-\ref{main:KLC}).

\begin{algorithm}[htp]
	\small
	\caption{Solving ACLP+RPDP}  \label{alg:main}
		
	\begin{algorithmic}[1]
			\State {\bf ACLP+RPDP} \,\, {\it in:}~$scenario, surplus, tmax, ntours$ \, \, {\it out:}~$answer$
			
			\State Let $M$ be the set of pallets ({\it cfr.} Table \ref{tab:larger}) \label{main:M} 
			
			\State Let $K$, $L$\/ and $C$ be according to the $scenario$ ({\it cfr.} Tables \ref{tab:costs}, \ref{tab:larger} and \ref{tab:scenarios}) \label{main:KLC} 
						
			\State $N \gets ItemsGeneration(scenario, surplus)$ \label{main:items}

			\For {each $method$}  \label{main:method} 
			
				\For {$i \gets 1$ to $ntours$} \label{main:loop1}

					\State $f_i \gets SolveTour(\pi[i], L, M, C, N, method, tmax/ntours)$ \label{main:fpi}
				
				\EndFor \label{main:loop2}
				
				\State $answer[scenario,surplus,method] \gets \max f_i$ \label{main:a2}
			\EndFor
	\end{algorithmic}

\end{algorithm}

$surplus$ is passed to $ItemsGeneration$\/ (line \ref{main:items}), responsible for creating the items to be shipped, which will be presented in the next section (Algorithm \ref{alg:itemsgen}). $tmax$\/ is the runtime limit, which will be distributed among the tours (line \ref{main:fpi}). $method$\/ corresponds to a MIP solver or a heuristic to the node-by-node solution $SolveTour$, which will be presented in subsection \ref{methods}. The best results obtained by testing all tours are stored in $answer$\/ (line \ref{main:a2}), which is the output of this algorithm.

$\pi[]$ is a vector of tours indexed from 1 to $ntours$. When $ntours=2$, $\pi[1]$\/ and $\pi[2]$\/ are the optimal solutions of the corresponding TSP. When $ntours=K!$, $\pi[i]$\/ is the $i^{th}$ permutation of $S_K$. Finally, when $ntours=100$, $\pi[1], \ldots, \pi[100]$\/ are the solutions provided by a TSP heuristic.

\begin{table}[h]
	\centering
	\caption{Testing scenarios}  \label{tab:scenarios}
	\begin{tabular}{c c c }
		\toprule
		{\bf Scenario} & $K$ & $L$  \\		
		\midrule
		1 & 2    & \{$0,1,2$\}            \\
		2 & 3    & \{$0,1,2,3$\}         \\
		3 & 4    & \{$0,1,2,3,4$\}        \\
		4 & 5    & \{$0,1,2,3,4,5$\}      \\
		5 & 6    & \{$0,1,2,3,4,5,6$\}    \\
		\bottomrule
	\end{tabular}
\end{table}

Next, we will present two subsections: in the first we explain how $SolveTour$\/ is executed. In the second we will present the heuristics developed for node-by-node solutions.

\subsection{SolveTour algorithm}
\label{tour}

As we commented in the previous subsection, we will adopt the strategy of presetting the destinations of each pallet throughout the tour. This is feasible in practical cases where $1<K<m$. For this, each pallet $i$\/ also has a field $T^k_i$, $0\leq k \leq K$, which stores its next destination after being loaded at node $k$. For this reason, $T^k_i \in L_k$, $1 \leq i \leq m$, $0\leq k \leq K$.

$SolveTour$\/ is described in Algorithm \ref{alg:tour}, where $\pi$\/ is a permutation of the nodes (excluding the base) that defines the order of visits in this tour, $method$\/ corresponds to a MIP solver or a heuristic for solving the node-by-node problems, and $tmax$\/ is the runtime limit of this tour.

\begin{algorithm}[H]
	\small
	\caption{Solving tour $\pi$ with $method$}  \label{alg:tour}
	\begin{algorithmic}[1]
		
		\State {\bf SolveTour} \,\, {\it in: $\pi, L, M, C, N, method$, $tmax$} \, \, {\it out: $score / cost$}
		
		\State $\pi_0     \gets 0$ \label{tour:pi1} \Comment{all tours start and end at the base}
		\State $\pi_{K+1} \gets 0$ \label{tour:pi2} 
		
		\State $score \gets 0$ \label{tour:score}
		\State $cost  \gets 0$ \label{tour:cost}

		\For {$k \gets 0$ to $K$} \label{tour:loop1}
		\State $t_{\pi_k} = (\sum_{j=1}^{n_k} v_j/ \sum_{k=0}^K \sum_{j=1}^{n_k} v_j) * tmax$ \label{tour:tnode}	 \Comment{runtime limit proportional to the shipment volume}
	
		\State $L_{\pi_k} \gets L - \{\pi_0,\pi_1,\ldots,\pi_k\}$  \label{tour:lk1}	 \Comment{the set of remaining nodes is updated}	
		\State $T_i^{\pi_k} \gets -1$, $1 \leq i \leq m$ \label{tour:-1}             \Comment{the pallet destination is unset}
		
		\If {$k = 0$}
			\State Let $G_1(M \cup N_0, \varnothing)$ \label{tour:g11}               \Comment{no packed contents at the base}
		\Else
			\State $E_{Q_{\pi_k}}, M \gets UpdatePacked(M, Q_{\pi_k}, \pi_k)$ \label{tour:dest}			
			\State Let $G_1(M \cup N_{\pi_k} \cup Q_{\pi_k}, E_{Q_{\pi_k}})$ \label{tour:g12}
		\EndIf  \label{tour:lk2}	
		
		\State $M \gets SetPalletsDestinations(M, \pi_k )$ \label{tour:dest2}
		
		\State $G_2 \gets SolveNode(method, \pi_k, G_1, t_{\pi_k})$  \label{tour:SolveNode}
		
		\State $s, \tau \gets ScoreAndTorque(\pi_k, G_2)$ \label{tour:analyse}
		
		\State $score \gets score + s$ \label{tour:score2}
		\State $cost  \gets cost  + c_{\pi_k,\pi_{k+1}} \times (1 + c_g \times |\tau|)$ \label{tour:cost2} 
		
		\EndFor  \label{tour:loop2}

	\end{algorithmic}
\end{algorithm}

As we mentioned in the previous section, all tours start and end at the base $0$\/ (lines \ref{tour:pi1}-\ref{tour:pi2}). After initializing the score and cost values (lines \ref{tour:score}-\ref{tour:cost}), there is a loop for the $K+1$\/ flights (lines \ref{tour:loop1}-\ref{tour:loop2}). Initially, the runtime limit $t_{\pi_k}$\/ for each node is calculated (line \ref{tour:tnode}), the set $L_{\pi_k}$\/ of remaining nodes is updated (line \ref{tour:lk1}), and the pallet destinations are unset (line \ref{tour:-1}).

When the aircraft is at the base, the initial graph $G_1$\/ is empty, and there are no packed contents (line \ref{tour:g11}). Otherwise, $UpdatePacked$\/ (line \ref{tour:dest}) returns the set of packed contents that have not yet reached their destination and remain on board, rearranging them on the pallets to minimize CG deviation. This allocation is stored in graph $G_1$\/ (line \ref{tour:g12}).

$SetPalletsDestinations$\/ (line \ref{tour:dest2}) presets the destination of each pallet based on the volume demands of the current node without changing the pallet's destination with packed contents.

Finally, $SolveNode$\/ includes the edges corresponding to the items shipped at the current node, returning the graph $G_2$\/ (line \ref{tour:dest2}). The score and the CG deviation of $G_2$\/ are calculated (line \ref{tour:analyse}) and accumulated (lines \ref{tour:score2}-\ref{tour:cost2}), allowing the final result of this tour as output.

$UpdatePacked$, described in Algorithm \ref{alg:cons}, finds the best packed-pallet allocation, in terms of CG deviation, for the packed contents that remain on board.

\begin{algorithm}[H]
	\small
	\caption{Updating the packed contents that remain boarded at node $\pi_k$}  \label{alg:cons}
	\begin{algorithmic}[1]
				
		\State {\bf UpdatePacked} \,\, {\it in: ${M, Q_{\pi_k}, \pi_k}$} \, \, {\it out: $E_{Q_{\pi_k}}, M$}
		
		\State $E_{Q_{\pi_k}} \gets MinCGDeviation(E_{Q_{\pi_k}})$   \label{cons:mincg}
		
		\For{$i \gets 1$ to $m$} \label{cons:Ybegin}
		
		\For{$q \gets 1$ to $m_{\pi_k}$}
		
		\State $T_i^{\pi_k} \gets -1$
		
		\If{$(i, q) \in E_{Q_{\pi_k}}$} 
		\State $T_i^{\pi_k} \gets to_q$   \label{cons:dest} \Comment{reassign pallet destinations}
		\EndIf
		\EndFor		
		\EndFor \label{cons:Yend}
				
	\end{algorithmic}
\end{algorithm}

$MinCGDeviation$\/ (line \ref{cons:mincg}) relocates the packed contents on the pallets, minimizing torque and ensuring that they all remain on board, one packed content on each pallet. It is run through a MIP solver with the objective function (\ref{eq:torque}) and the constraints (\ref{eq:embarked}) and (\ref{eq:one}). As there are few variables, $E_{Q_{\pi_k}}$\/ is obtained in less than $30$\/ milliseconds. Finally, the destination of each pallet with packed content is updated (lines \ref{cons:Ybegin}-\ref{cons:Yend}).

\begin{equation} \label{eq:torque}
	\min f =  \Big | \sum_{i=1}^{m} \sum_{q=1}^{m_{\pi_k}} Y^k_{iq} \times w_q \times D_i^{long}  \Big |
\end{equation}

\begin{equation} \label{eq:embarked}
	\sum_{i=1}^{m} Y^k_{iq} = 1;\ q \in \{1,\ldots,m_{\pi_k}\}
\end{equation}

\begin{equation} \label{eq:one}
	\sum_{q=1}^{m_{\pi_k}} Y^k_{iq} \leq 1;\ i \in \{1,\ldots,m\}
\end{equation}

$SetPalletsDestinations$, which sets the pallets destination not yet defined, is described in Algorithm \ref{alg:dest}.

\begin{algorithm}[H]
	\small
	\caption{Setting pallets destination based on the items to be embarked at node $\pi_k$}  \label{alg:dest}
	\begin{algorithmic}[1]
		
		\State {\bf SetPalletsDestinations} \,\, {\it in: $M, \pi_k$} \, \, {\it out: $M$}
		\State $vol_x \gets  0$, $x \in L_{\pi_k}$  \label{dest:vector1}
		\State $max \gets 0$  \Comment{destination with maximum volume demand}
		\State $total \gets 0$ 
		\For{$j \gets 1$ to $n_{\pi_k}$}
		\If {$to_j \in L_{\pi_k}$} 
		\State $vol_{to_j} \gets vol_{to_j} + v_j$ 
		\State $total \gets total + v_j$ 
		\If {$vol_{to_j} > vol_{max}$}
		\State $max \gets to_j$ \label{dest:max} 
		\EndIf
		\EndIf
		\EndFor
		\For{$x \in L_{\pi_k}$} \label{dest:propor1}
		\If {$vol_{x} \neq 0$}
		\State $needed \gets \max \{1, \lfloor{ (m-m_{\pi_k}) \times vol_{x}/total}\rfloor \}$  
		\label{dest:needed}
		\State $np \gets 0$
		\For{$i \gets 1$ to $m$}
		\If{($np < needed$) {\bf and} ($T_i^{\pi_k} = -1$)}  
		\State $T_i^{\pi_k} \gets x$
		\State $np \gets np + 1$ \Comment{number of necessary pallets to node $x$}
		\EndIf
		\EndFor
		\EndIf 
		\EndFor \label{dest:propor2}
		\For{$i \gets 1$ to $m$} \label{dest:max1}
		\If{$T_i^{\pi_k} \gets -1$} 
		\State $T_i^{\pi_k} \gets max$  \Comment{any remaining pallet is assigned to the maximum demand destination}
		\EndIf
		\EndFor \label{dest:max2}
		
	\end{algorithmic}
\end{algorithm}

$vol$\/ stores the demand volume of items destined for the non-visited nodes (line \ref{dest:vector1}). The destination of empty pallets is defined proportionally to the volume of items to be embarked (lines \ref{dest:propor1}-\ref{dest:propor2}). $max$\/ is the destination with maximum volume demand (line \ref{dest:max}), and $needed$\/ is the number of necessary pallets to node $x$ (line \ref{dest:needed}). The destination with the maximum volume defines any remaining pallets (lines \ref{dest:max1}-\ref{dest:max2}).

$ScoreAndTorque$, described in Algorithm \ref{alg:eval}, evaluates the allocation graph $G$\/ generated by $SolveNode$\/ at node $\pi_k$\/ and returns the corresponding cargo score and aircraft torque.

\begin{algorithm}[H]
	\small
	\caption{Cargo score and aircraft torque}  \label{alg:eval}
	\begin{algorithmic}[1]
				
		\State {\bf ScoreAndTorque} \,\, {\it in: $\pi_k, G$} \, \, {\it out: $s, \tau$}
		
		\State Let $G(V_{\pi_k}, E_{Q_{\pi_k}} \cup E_{N_{\pi_k}})$
		
		\State $s \gets 0$
		\State $\tau_i \gets 0$, $1 \leq i \leq m$
		\For{$i \gets 1$ to $m$} \label{eval:loop1}
		
		\For{$j \gets 1$ to $n_{\pi_k}$}
		
		\If{$X_{ij}^{\pi_k} = 1$}
		\State $s \gets s + s_j$ \label{eval:score1} \Comment{accumulates cargo score}
		\State $\tau_i \gets \tau_i + w_j \times D_i^{long}$  \label{eval:eps1} \Comment{accumulates aircraft torque}
		\EndIf
		\EndFor	
		
		\For{$q \gets 1$ to $m_{\pi_k}$}
		\If{$Y_{iq}^{\pi_k} = 1$}
		\State $s \gets s + s_q$   \label{eval:score2} \Comment{accumulates cargo score}
		\State $\tau_i \gets \tau_i + w_q \times D_i^{long}$ \label{eval:eps2} \Comment{accumulates aircraft torque}
		\EndIf
		\EndFor	   
		\EndFor 		 \label{eval:loop2}
		\State $\tau \gets  \sum_{i=1}^{m} \tau_i / (W_{max} \times limit^{CG}_{long})$ \label{eval:torque} \Comment{final calculation of the aircraft torque}
				
	\end{algorithmic}
\end{algorithm}

Algorithm \ref{alg:eval}\/ consists of a loop that goes through all the pallets (lines \ref{eval:loop1}-\ref{eval:loop2}), accumulating the scores (lines \ref{eval:score1} and \ref{eval:score2}) and the torques (lines \ref{eval:eps1} and \ref{eval:eps2}) of the shipped items, allowing the final calculation of the aircraft torque (line \ref{eval:torque}).

\subsection{Node-by-node solutions}
\label{methods}

In this subsection, we present two implementations of $SolveNode$\/ algorithm: with a MIP solver and with heuristics.

\subsubsection{Node-by-node solutions with a MIP solver}
\label{solver}

Our strategy adopted in $SolveTour$\/ defines the values of some variables: the set of nodes to be visited is updated, the packed contents that remain on board are reallocated to minimize the CG deviation, and the pallet's destinations are determined according to the volume of items available for shipment. 

In this way, the mathematical model for $SolveNode(MIP,\pi_k,G, tmax)$\/ becomes simpler, which finds an allocation of available items at node $\pi_k$\/ using previously defined values of $L_{\pi_k}$, $T_i^{\pi_k}$, and $a^{\pi_k}_q$. Thus, we use a MIP solver with a runtime limit $tmax$\/ at node $\pi_k$ to maximize the objective function (\ref{eq:obj2}) with the calculus equations (\ref{eq:newf}) to (\ref{eq:costs2}), subject to the constraints (\ref{eq:tau2}) to (\ref{eq:ifB}). The binary variables $X_{ij}$\/ and $Y_{iq}$\/ define the sets of edges $E_{N_{\pi_k}}$\/ and $E_{Q_{\pi_k}}$, respectively, included in graph $G$.

\begin{equation} \label{eq:obj2}
	\max f= \tilde{s} / \tilde{c}
\end{equation}

\begin{equation} \label{eq:newf}
	\tilde{s} = \sum_{i=1}^{m} \sum_{j=1}^{n_{\pi_k}} X_{ij} \times s_j
\end{equation}

\begin{equation} 
	\tau_{\pi_k} = \sum_{i=1}^{m} \Big [ D_i^{long} \times (\sum_{j=1}^{n_{\pi_k}} X_{ij} \times w_j +  \sum_{q=1}^{m_{\pi_k}} Y_{iq} \times w_q) \Big ] \Big / W_{max} \times limit^{CG}_{long}
\end{equation}

\begin{equation} \label{eq:costs2}
	\tilde{c} =  c_{\pi_k, \pi_{k+1}}\times(1+c_g\times|\tau_{\pi_k}|)
\end{equation}

\begin{equation} \label{eq:tau2}
	|\tau_{\pi_k}| \leq 1
\end{equation}

\begin{equation} 
	\sum_{j=1}^{n_{\pi_k}} X_{ij} \times w_j + \sum_{q=1}^{m_{\pi_k}} Y_{iq} \times w_q  \leq W_i; \ i \in \{1, \ldots, m\}
\end{equation}

\begin{equation} 
	\sum_{j=1}^{n_{\pi_k}} X_{ij} \times v_j + \sum_{q=1}^{m_{\pi_k}} Y_{iq} \times v_q  \leq\ V_i; \ i \in \{1, \ldots, m\}
\end{equation}

\begin{equation} 
	\sum_{i=1}^{m} X_{ij} \leq 1; \ j \in \{1, \ldots, n_{\pi_k}\}
\end{equation}

\begin{equation} \label{eq:if2}
	X_{ij} = 0;\ to_j \notin L_{\pi_k}; \ i \in \{1, \ldots, m\}; \ j \in \{1, \ldots, n_{\pi_k}\}
\end{equation}

\begin{equation} \label{eq:ifA}
	X_{ij} \leq X_{ij} \times (T_i^{\pi_k} - to_j + 1); \ i \in \{1, \ldots, m\}; \ j \in \{1, \ldots, n_{\pi_k}\}
\end{equation}

\begin{equation} \label{eq:ifB}
	X_{ij} \leq X_{ij} \times (to_j - T_i^{\pi_k} + 1 ); \ i \in \{1, \ldots, m\}; \ j \in \{1, \ldots, n_{\pi_k}\}
\end{equation}

The constraints (\ref{eq:ifA}) and (\ref{eq:ifB}) are equivalent to $X_{ij} =1$\/ if $to_j = T_i^{\pi_k}$, and $X_{ij} =0$\/ otherwise .

\subsubsection{Node-by-node solutions with heuristics}

One of the main objectives of this work was to find a quick heuristic that offers a good-quality solution for the node-by-node problem. Taking this into account, we design algorithms based on five known meta-heuristics: {\it Ant Colony Optimization} (ACO) \citep{Dorigo1992,DorigoManiezzoColorni1996}, {\it Noising Method Optimization} (NMO) \citep{CharonHudry1993,CharonHudry2001,Zhan2020}, {\it Tabu Search} (TS) \citep{Glover1986}, {\it Greedy Randomized Adaptive Search Procedure} (GRASP) \citep{FeoResende1989}, and {\it Genetic Algorithm} (GA) \citep{holland1992adaptation}. We considered several ideas from the literature \citep{NiarFreville1997, Fidanova2006, Alonso2019, Zhan2020, PeerlinckAmy2022MFEO}, and we were careful to use the same data structures and procedures in all implementations to enforce fair results comparison.

However, the heuristic that presented better solutions was none of the previous ones. In this subsection, we will present a new heuristic for the node-by-node problem, called {\it Shims}. Like in mechanics, shims are collections of spacers to fill gaps, which may be composed of parts with different thicknesses (see Figure \ref{fig:shims}). This strategy is based on a practical observation: usually, subsets of smaller and lighter items are saved for later adjustments to the remaining available space.

The selection of edges for $E_{ N_{\pi_k} }$ uses the {\it edge attractiveness}\/ $\theta_{ij}$\/ (\ref{eq:edge}), which can be understood as the tendency to allocate item $j$\/ to pallet $i$\/ at node $\pi_k$. It is directly proportional to the score, and inversely to the volume and the torque of each item.

\begin{equation} \label{eq:edge}
	\theta_{ij} = \frac{ s_j }{ v_j }\ \times \Big ( 1 - \frac{ w_j \times |D_i^{long}| }{ \max w_j \times\  \max | D_i^{long} | } \Big );\ i \in \{1,\ldots,m\},\ j \in \{1,\ldots,n_{\pi_k}\}
\end{equation}

\begin{table}[H]
	
	\begin{minipage}{0.08\linewidth}
		
	\end{minipage}\hfill
	\begin{minipage}{0.34\linewidth}
		
		\includegraphics[scale=0.25]{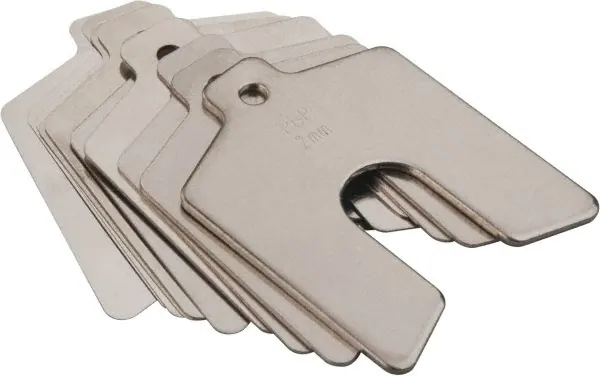}
		\captionof{figure}{{\it Shims} of various thicknesses}
		\tiny\textsuperscript{Source: www.mscdirect.com/product/details/70475967}
		\label{fig:shims}
		
	\end{minipage}\hfill
	\begin{minipage}{0.58\linewidth}
		\scriptsize
		\begin{tikzpicture}[scale=0.75, samples=100]
			\filldraw[fill=blue!3!white, draw=black] (0, 0) rectangle (12.5, 5.5);
			\draw[->] (.2, .6) -- coordinate (x axis mid) (12, .6);
			\node at (5, 0.3) {$\eta_1$};
			\node at (5, 2.5) {|};
			\node at (6, 2.5) {\small{shims}};
			\node at (7, 2.5) {|};
			\node at (7, 0.3) {$\eta_2$};
			\node at (0.3, 5) {$\theta_{ij}$};
			\draw[->] (.6, .2) -- coordinate (y axis mid) (0.6, 5.3);
			\node at (12, 0.3) {$e_{ij}$};
			\draw[smooth, domain = 0.09:2, color=black, thick] plot (.3+1/\x,{4.2+log2(\x)});
			\draw[->, dashed] (5, 0.6)--(5, 2.0);
			\draw[->, dashed] (7, 0.6)--(7, 1.5);
		\end{tikzpicture}
		
		
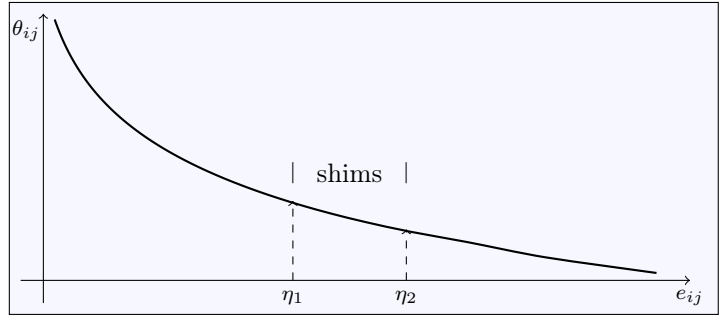
\captionof{figure}{$n_{\pi_k}$\/ possible edges $e_{ij}$\/ sorted by $\theta_{ij}$\/ in non-ascending order}
		\label{fig:whip}		
	\end{minipage}
\end{table}

Considering only the items that can be shipped at node $\pi_k$, Figure \ref{fig:whip} represents $n_{\pi_k}$\/ possible edges $e_{ij}$\/ of pallet $i$\/ sorted by $\theta_{ij}$\/ in non-ascending order. Initially, {\it Shims}\/ builds a greedy solution for pallet $i$\/ selecting edges up to index $\eta_1$ ({\it greedy phase}). Then, with the edges between $\eta_1$\/ and $\eta_2$, it elaborates different possible complements ({\it composition phase}), including later the best ones in the same pallet ({\it selection phase}). {\it Shims}\/ is depicted in Algorithm \ref{alg:shims}.

\begin{algorithm}[H]
\small
	\caption{{\it Shims}\/ heuristic at node $\pi_k$}  \label{alg:shims}
	\begin{algorithmic}[1]
				
		\State {\bf SolveNode} \,\, {\it in: {\it Shims}, $\pi_k, G, tmax, level_1, level_2$} \, \, {\it out: $G(M \cup N_{\pi_k} \cup Q_{\pi_k}, E_{Q_{\pi_k}} \cup E_{N_{\pi_k}})$}
		
		\State $T_{begin} \gets$ current system time
		
		\State Let $G(M \cup N_{\pi_k} \cup Q_{\pi_k}, E_{Q_{\pi_k}})$   \label{shims:initQ}
		
		\State Sort $M$ by $|D_i^{long}|$ in non-descending order \label{shims:pallets}
		
		\State $E_{N_{\pi_k}} \gets \varnothing$ 
		
		\State $\tau_{max} \gets W_{max} \times limit^{CG}_{long}$

		\For{$i \gets 1$ to $m$}  \label{shims:eta1_a}
		
		\State $\tau_{\pi_k} \gets \sum_{(i,q)\in E_{Q_{\pi_k}}} w_q \times D_i^{long}$ 
		\State $vol_i \gets \sum_{(i,q)\in E_{Q_{\pi_k}}} v_q$ 
		
		\State Let $E$\/ be an array of $n_{\pi_k}$ possibles edges of pallet $i$\/ sorted by $\theta_{ij}$\/ in non-ascending order \label{shims:possible}
		\State $\eta_1 \gets 1$  
		\Repeat
		\State $e_{ij} \gets E_{\eta_1}$ 
		
		\If{ ($E_{N_{\pi_k}} \cup \{e_{ij}\}$ is feasible) {\bf and} ($vol_i \leq V_i \times level_1$) {\bf and} ($| \tau_{\pi_k} + w_j \times D_i^{long} | \leq W_{max} \times limit^{CG}_{long}$)} 
		
		\State $E_{N_{\pi_k}} \gets E_{N_{\pi_k}} \cup \{e_{ij}\}$
		
		\State $vol_i \gets vol_i + v_j$
		
		\State $\tau_{\pi_k} \gets \tau_{\pi_k} + w_j \times D_i^{long}$
		
		\State $\eta_1 \gets \eta_1 + 1$ 
		
		\EndIf
		\Until ($vol_i > V_i \times level_1$) {\bf or} ($\eta_1>n_{\pi_k}$) \label{shims:endfirst}	
		
		\State $slack_i \gets V_i - vol_i$  \label{shims:beginsecond}
		\State $\eta_2 \gets \eta_1$  \label{shims:eta2a}
		\While{($\eta_2 \leq n_{\pi_k}$) {\bf and} ($ vol_i < V_i \times level_2$)}
	
		\State $e_{ij} \gets E_{\eta_2}$
		
		\State $vol_i \gets vol_i + v_j$
		\State $\eta_2 \gets \eta_2 + 1$  \label{shims:eta2} 
		\EndWhile \label{shims:endsecond} 

		\State $vol \gets 0$; $b \gets 1$; $shims_b \gets \varnothing$; $Set \gets \{ shims_b \}$  \label{shims:beginthird} 
				
		\For{$x \gets \eta_1$ to $\eta_2$} \label{edges_indexes}
		
		\If{ $T_{current} - T_{begin} > tmax$}
		\State {\bf break}
		\EndIf
		
		\State $NewShims \gets$ {\bf True} \label{new_shims}
		\State $e_{ij} \gets E_x$
		
		\For{$shims \in Set$} \label{shims_set}
		
		\If {($e_{ij} \not\in (E_{N_{\pi_k}} \cup shims))$ {\bf and} ($e_{ij}$ is feasible)  {\bf and} ($(v_j + vol) \leq slack_i$)}
		
		\State $shims \gets shims \cup \{e_{ij}\}$
		\State $vol \gets vol + v_j$
		\State $NewShims \gets$ {\bf False} \label{new_shims_false}
		
		\State {\bf break}  
		\EndIf
		
		\EndFor 
		
		\If{$NewShims$} \label{new_shims2}
		\State $vol \gets 0$; $b \gets b + 1$;  $shims_b \gets \{e_{ij}\}$
		\State $Set \gets Set \cup \{ shims_b \}$
		\EndIf
		
		\EndFor 
		
		\State $sh_w \gets shims$, where $shims \in Set$ and $\sum_{e_{ij} \in shims} w_j$\/ is maximum  \label{best_weight}		
		\State $sh_v \gets shims$, where $shims \in Set$ and $\sum_{e_{ij} \in shims} v_j$\/ is maximum \label{best_volume}
		\State $sh_{best} \gets shims$, where $shims \in \{sh_w, sh_v\}$\/ and $\sum_{e_{ij} \in shims} s_j$\/ is maximum \label{best_score}
		\State $E_{N_{\pi_k}} \gets E_{N_{\pi_k}} \cup sh_{best}$ \label{shims:endthird}  
		
		\EndFor

	\end{algorithmic}
\end{algorithm}

$tmax$\/ is the runtime limit for {\it Shims}, and $level_1$\/ and $level_2$\/ are volume thresholds for indices $\eta_1$ and $\eta_2$, respectively. Initially, $Q_{\pi_k}$\/ (line \ref{shims:initQ}) corresponds to the packed contents that remain on board. It is important to remember that $E_{Q_{\pi_k}}$\/ and $M$\/ were modified by the procedure $UpdatePacked(M, Q_{\pi_k}, \pi_k)$ and the procedure $SetPalletsDestinations(M, \pi_k)$. Then, the pallets are considered in non-descending order of $|D_i^{long}|$.

For each pallet $i$, its $n_{\pi_k}$\/ possible edges $e_{ij}$\/ are considered in non-increasing order of $\theta_{ij}$:

\begin{itemize}

	\item In the {\it greedy phase}\/ (lines \ref{shims:pallets}-\ref{shims:endfirst}), a partial solution for each pallet $i$\/ is constructed by adding edges following $\theta_{ij}$\/ order. Indices $\eta_1$\/ and $\eta_2$\/ refer to the accumulated volumes $V_i \times level_1$\/ and $V_i \times level_2$, respectively, which were defined empirically by the {\it irace}\/ tool \citep{LopezIbanezManuel2016}. We will explain this in subsection \ref{items}.
	
	\item In the {\it composition phase}\/ (lines \ref{shims:beginsecond}-\ref{shims:endsecond}), a set of shims named $Set$\/ is created for each pallet $i$, where each shim is formed by a set of edges in the range $[\eta_1,\eta_2]$, whose total volume is limited by $slack_i$. In this phase, the heuristic that provided the best results, both in terms of time and quality, is based on {\it First-Fit Decreasing}, which is an approximation algorithm for the {\it Bin Packing Problem}\/ \citep{JohnsonGarey1985}. Basically, shims are created by accumulating the following edges, taking $slack_i$\/ as a limit.
	
	\item In the {\it selection phase}\/ (lines \ref{shims:beginthird}-\ref{shims:endthird}), the best shim in $Set$\/ is chosen. Initially, two shims are found: $sh_w$\/ with larger weight and $sh_v$\/ with larger volume. Between the two, the one with the highest score will be chosen, and its edges will be inserted into $E_{N_{\pi_k}}$.
	
\end{itemize}

\subsection{Time complexity using {\it Shims}}

We finish this section with the analysis of the time complexity of our ACLP+RPDP solution, considering the use of the {\it Shims}\/ heuristic.

In this process, when $K>6$\/ (that is, in unusual cases of air transportation), we need a GA-based TSP heuristic to generate 100 tours of size $K$. This heuristic can be chosen and calibrated to be fast enough: specifically, in the tests we will present in the next section, it took just $33 s$ with $K=15$. For this reason, we will not perform its complexity analysis. We will also not consider the time spent by the {\it irace}\/ tool, in defining parameters for {\it Shims}, as it is only executed once in the calibration of our method.

Let $n = \max n_k$. Without loss of generality, we assume that $m = \mathcal{O}(n)$ and therefore $\max m_k = \mathcal{O}(n)$. This way, to read the input data, Algorithm \ref{alg:main} requires time $\mathcal{O}(m + K^2 + K.n + m.\max m_k) = \mathcal{O}(m.n)$, because $K<m$. 

Since $ntours \leq 6!$, it is enough to calculate the complexity of Algorithm \ref{alg:tour}:

\begin{itemize}
	
	\item Its initial variables are $\mathcal{O}(K.n)$.
	
	\item Algorithm \ref{alg:cons} minimizes the CG deviation related to packed contents using a MIP solver. We are unable to analyse its time complexity, but we found that it is very fast in all cases considered in our work, where $m=18$: it spent a maximum of $30 ms$. The final loop of Algorithm \ref{alg:cons} takes time $\mathcal{O}(m.n)$.
	
	\item Algorithm \ref{alg:dest} takes $\mathcal{O}(n+K+m)$ time.
	
	\item {\it Shims}, described in Algorithm \ref{alg:shims}, spends $\mathcal{O}(m.\log m)$ in ordering the pallets and then performs a loop with $m$\/ iterations. Each of these iterations spends $\mathcal{O}(n.\log n)$\/ on ordering the edges by the value of $\theta_{ij}$, and $\mathcal{O}(n)$\/ on the others statements. So its total time is $\mathcal{O}(m.n.\log n)$.
	
	\item Algorithm \ref{alg:eval} takes $\mathcal{O}(m.n)$ time.
	
\end{itemize}

Therefore, since there are $K$\/ iterations in Algorithm \ref{alg:tour}, we can conclude that its total time is $\mathcal{O}(K.m.n.\log n)$, which is the time complexity of  ACLP+RPDP solution.

\section{Implementation and results}
\label{results}

This section is composed of two parts: the generation of the test instances and the results obtained in our implementation.

\subsection{Instances generation}
\label{items}

As we are dealing with a new problem that until now had not been modelled in the literature, we have to create our own benchmarks. For this, we based it on the characteristics of real airlifts carried out by the {\it Brazilian Air Force}, as described below.

In the delivery of supplies carried out in Brazil from 2008 to 2010, 23\% of the items weighed between $10kg$ and $20kg$, 22\% from $21kg$ to $40kg$, 24\% from $41kg$ to $80kg$, 23\% from $81kg$ to $200kg$, and 8\% between $201kg$ and $340kg$. These five groups of items are described in Table \ref{tab:weights}, where $P$\/ represents the group probability. On the other hand, the average density of these items is approximately $246 kg/m^3$.

\begin{table}[h]
	\centering
	\caption{Items weight distribution}  \label{tab:weights}
	\begin{tabular}{c c c c }
		\toprule
		$item$ & $P$ & $low$ ($kg$) & $high$ ($kg$) \\		
		\midrule
		1              & 0.23           & 10  & 20  \\
		2              & 0.22           & 21  & 40  \\
		3              & 0.24           & 41  & 80  \\		
		4              & 0.23           & 81  & 200 \\
		5              & 0.08           & 201 & 340 \\
		\bottomrule
	\end{tabular}
\end{table}

In the generation of test instances, we use two types of random selections:

\begin{itemize}

	\item $RandomInt(i_1,i_2)$: randomly selects a integer number in $[i_1,i_2]$, where $i_1$ and $i_2$ are integer numbers;
	
	\item $Roulette(set)$ biased through $\phi$: selects an element from $set$, where the probability of each element is proportional to the value of a given function $\phi$\/ defined on $set$.
	
\end{itemize}

The procedure $ItemsGeneration$, which generates $N$ (all items to be moved among the nodes), is described in Algorithm \ref{alg:itemsgen}.

\begin{algorithm}[H]
	\small
	\caption{Generating items}  \label{alg:itemsgen}
	\begin{algorithmic}[1]
		\State {\bf ItemsGeneration} \,\, {\it in: $scenario, surplus$} \, \, {\it out: $N$}
		\State Let $L$ be the set of nodes and $M$ the set of pallets  \label{ig:LM}
		\State $limit \gets surplus \times \sum_{i=1}^{m} V_i$ \label{ig:extended}
		\For{$k \gets 0$ to $K$}
		\State $N_k \gets \varnothing$
		\State $j \gets 0$
		\State $vol \gets 0$	\label{ig:totals}	
		\While{$vol < limit$}
		\State $j \gets j+1$
		\State Let $t_j^k$\/ be the item $j$\/ at the node $k$
		\Repeat
		\State $to_j \gets RandomInt(0, K)$ \label{ig:dest}
		\Until{$to_j \neq k$}
		\State $x = Roulette(item)$ biased through $P$ \label{ig:weight1} \Comment{From Table \ref{tab:weights}}	 
		\State $w_j \gets RandomInt(low(x), high(x))$        \label{ig:weight2}		
		\State $s_j \gets \round{100 \times (1 - \log_{10}(RandomInt(1, 9)))} $ \label{ig:score}
		\State $v_j \gets w_j / RandomInt(148, 344)$ \label{ig:volume}
		\State $vol \gets vol + v_j$ 
		\State $N_k \gets N_k \cup \{t_j^k\}$ 
		\EndWhile
		\EndFor
		\State $N \gets \bigcup_{0 \leq k \leq K} N_k$
	\end{algorithmic}
\end{algorithm}

$scenario$\/ defines $L$\/ and $M$\/ (line \ref{ig:LM}), and the argument $surplus$\/ sets a limit on the total volume of items at each node (line \ref{ig:extended}). To avoid simply loading all items, we use $surplus \in \{1.2,\ 1.5,\ 2.0\}$\/. This also represents more instances for tests in each scenario.

For each generated $t^k_j$ item, its destination is randomly selected (line \ref{ig:dest}), its weight has a distribution according to Table \ref{tab:weights} (lines \ref{ig:weight1}-\ref{ig:weight2}), its score varies $100$\/ (highest) and $5$\/ (lowest) according to a logarithmic scale (line \ref{ig:score}), and its volume is randomly defined from the density, where we allow a variation of 40\% around the average density of $246 kg/m^3$\/ (line \ref{ig:volume}).

To determine the parameters $level_1$\/ and $level_2$\/ used by \textit{Shims}, we previously carried out some experiments with the {\it irace}\/ tool \citep{LopezIbanezManuel2016}, the results of which are presented in Table \ref {tab:irace}. In these tests, of every 7 instances generated for each value of $surplus$, 4 were used as the training set and 3 as the testing set. We provided the ranges [0.8, 1.0] and [1.0, 2.0] for $level_1$\/ and $level_2$, respectively. In each experiment, there was a maximum of 3,000 runs so that \textit{irace}\/ would have enough data for its statistical tests. For more details, see {\tt cran.r-project.org/web/packages/irace/}.

	\begin{table}[H]
	\centering
	\caption{\textit{irace}\/ results}  \label{tab:irace}
	\begin{tabular}{cccc}
		\toprule
		
		$surplus$ & $level_1$ & $level_2$ & runtime (min)\\

		\midrule
		1.2& 0.8621 & 1.0539 & 47 \\
		1.5& 0.9199 & 1.1399 & 59 \\
		2.0& 0.9617 & 1.5706 & 63 \\						
		\bottomrule
	\end{tabular}
	\normalsize
\end{table}

\subsection{Obtained solutions: quality and runtimes}

In the tests performed, we used a 64-bit, 16 GB, 3.6 GHz, eight-core processor with {\it Linux Ubuntu 22.04.1 LTS 64-bit} as the operational system and {\it Python 3.10.4}\/ as the programming language. We also used the well-known solver {\it Gurobi}\/ ({\tt www.gurobi.com}), version 9.5.2.

We will first present the results obtained in practical cases, when $1<K\leq6$. Next, we will show how {\it Shims}\/ remains robust when $K>6$.

\subsubsection{Results when $1<K\leq6$}

We ran Algorithm \ref{alg:main} considering the 5 scenarios from Table \ref{tab:scenarios}, 3 values for  $surplus$\/ from $\{1.2, 1.5, 2.0\}$, 4 values for $tmax$\/ from $\{240s, 1200s, 2400s, 3600s\}$, and 7 different methods for the node-by-node solution: {\it Gurobi}\/ (subsection \ref{solver}), ACO, NMO, TS, GRASP, GA, and {\it Shims}\/ (Algorithm \ref{alg:shims}).

For {\it Gurobi}\/ to be able to solve the largest possible number of tests without memory overflow, we set its parameter {\it MIPgap}\/ to 1\%. This shortens its runtime, in addition to ensuring that its objective function $f$\/ is at most 1\% of the optimal solution. For more details, see {\tt www.support.gurobi.com}. For each $scenario$, $surplus$\/ and $tmax$\/ tested, 7 different instances were generated. 

Table \ref{tab:overall} succinctly shows the overall performance of the methods for the node-by-node solution. Only {\it Shims}\/ found a solution for all scenarios, as well as being the fastest.

\begin{table}[H]
	
	\centering
	\caption{Overall results}  \label{tab:overall}
	\small
	\begin{tabular}{cccc}	
		\toprule
		{\bf Method} & {\bf Best scenarios} & {\bf Worst scenarios} & {\bf Worst runtimes (min)}  \\
		\midrule
		NMO           & 4       & 5    & 60  \\
		ACO           & 2, 3    & 4, 5 & 25, 61  \\
		GRASP         & 1       & 4, 5 & 28, 55  \\
		TS            & -       & 5    & 44  \\
		GA            & -       & 1, 2, 3, 4, 5  & did not solve  \\
		\midrule
		{\it Gurobi}  & 1, 2, 3, 4    & 5    & did not solve  \\
		{\it Shims}   & 1, 2, 3, 4, 5 & -    & 3.26 \\
		\bottomrule
	\end{tabular}
	\normalsize
	
\end{table}

Table \ref{tab:mh} shows a particular case ($surplus=1.2$, $ntours=K!$, and $tmax=3600s$) in which the methods can find solutions for all scenarios. As can be seen, {\it Shims}\/ always obtained the best value for the objective function, in addition to being the fastest.

All methods used (ACO, NMO, TS, GRASP, and GA) generate a large number of solutions that require further evaluation, resulting in longer runtimes. On the other hand, {\it Shims}\/ is a constructive heuristic that continually builds a feasible solution, which makes it much faster. {\it Shims}\/ uses the {\it First-Fit Decreasing}\/ algorithm and the {\it irace}\/ tool, which follow a greedy process in the search for solution quality.

It is important to highlight that the speed of obtaining a balanced allocation at a node is essential for the ACLP+RDPD solution. Only in this way will it be possible to make the various comparisons between different routes, allowing the obtaining of an efficient flight itinerary with pickup and delivery.

\begin{table}[H]
	\centering
	\caption{Solutions with $surplus=1.2$, $ntours=K!$, and $tmax=3600s$}  \label{tab:mh}
	\setlength{\tabcolsep}{4.5pt}
	\footnotesize
	\begin{tabular}{c| cc|cc | cc| cc  | cc | cc }
		\toprule
		& \multicolumn{2}{c}{ NMO}     & \multicolumn{2}{c}{ ACO}          & \multicolumn{2}{c}{GRASP}     & \multicolumn{2}{c}{ TS} & \multicolumn{2}{c}{ GA} & \multicolumn{2}{c}{\it Shims}\\
		$scenario$ & $f$ & time (s)     & $f$ & time (s)          & $f$ & time (s)       & $f$ & time (s) & $f$ & time (s) & $f$ & time (s)\\	
		\midrule	
		1               &      8.38   &     4     &     8.40  &   17      &   8.40    &     12    &    8.39   &   4       &    5.03	  &    154    & 8.49& 1 \\		
		2               &     11.21   &     17    &    11.29  &   56      &   11.24   &     47    &   11.12   &   18      &    5.66	  &    544    &12.10& 2 \\		
		3               &     13.12   &     87    &    13.25  &   278     &   13.13   &     258   &   13.03   &   90      &    5.63	  &    2,614   &13.30& 8 \\
		4               &     13.32   &     520   &    13.73  &   1,516    &   13.31   &     1,690  &   13.22   &   586     &    5.45	  &    2,924   &14.49& 10 \\
        5               &     52.20   &     3,582  &    52.36  &   3,602    &   52.16   &     3,292  &   51.21   &   2,627    &    18.74  &    3,192   &53.61& 36 \\ 
					
		\bottomrule		
	\end{tabular}
	\normalsize
\end{table}

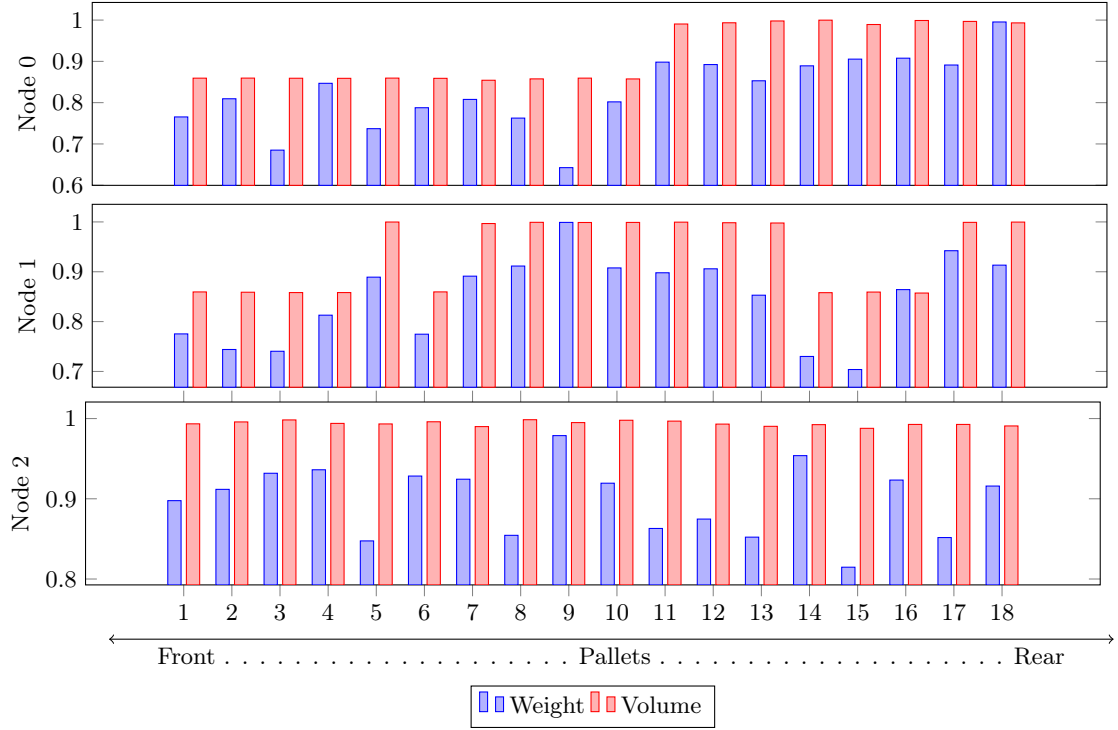
\begin{figure}[htp]
	\centering	
	\footnotesize
	
	\pgfplotsset{width=15cm,height=4cm,compat=1.13}
	\begin{tikzpicture}
		\begin{axis}[
			xtick=\empty,
			ybar,
			enlargelimits=0.12,
			ylabel={Node 0},
			symbolic x coords={1,2,3,4,5,6,7,8,9,10,11,12,13,14,15,16,17,18},
			bar width=5pt,
			]
			
			\addplot coordinates {(1,0.7654138319766367)(2,0.8094308769437574)(3,0.6850199132095408)(4,0.8468463660420169)(5,0.7369893385188176)(6,0.7877157566954458)(7,0.8078374028365509)(8,0.7625929630506392)(9,0.6425980437410217)(10,0.8019085281114264)(11,0.8981144246398084)(12,0.8922746169744865)(13,0.8529681500451808)(14,0.8891665195473366)(15,0.9053959007861928)(16,0.9076987362194541)(17,0.8911561631588362)(18,0.9953993552024535)};

			\addplot coordinates {(1,0.8593243243243244)(2,0.8594594594594593)(3,0.8591216216216214)(4,0.8589189189189187)(5,0.8594594594594593)(6,0.8589189189189191)(7,0.8542567567567567)(8,0.8575675675675676)(9,0.8593918918918919)(10,0.8574324324324325)(11,0.9904729729729729)(12,0.9934459459459459)(13,0.9978378378378379)(14,0.9997972972972973)(15,0.9892567567567567)(16,0.9989864864864866)(17,0.9967567567567567)(18,0.9931756756756758)};	
			
		\end{axis}
	\end{tikzpicture}
	
	\begin{tikzpicture}
		\begin{axis}[
			xtick=\empty,
			ybar,
			enlargelimits=0.12,
			legend style={at={(0.5,-0.2)}, anchor=north,legend columns=-1}, 
			ylabel={Node 1},
			symbolic x coords={1,2,3,4,5,6,7,8,9,10,11,12,13,14,15,16,17,18},
			bar width=5pt,
			]
			
			\addplot coordinates {(1,0.7753115232339638)(2,0.7440372483226444)(3,0.7404687619589224)(4,0.8128810527891527)(5,0.8891665195473366)(6,0.7747138784769633)(7,0.8911561631588362)(8,0.9114634506632016)(9,0.9990350016040416)(10,0.9076987362194541)(11,0.8979761771192329)(12,0.9059251895751423)(13,0.8529681500451808)(14,0.7300589473230157)(15,0.7038498441692468)(16,0.8641814862404096)(17,0.9422229510226497)(18,0.9132675161355484)};

			\addplot coordinates {(1,0.8595270270270271)(2,0.8589189189189188)(3,0.8582432432432435)(4,0.8582432432432434)(5,0.9997972972972973)(6,0.8595945945945944)(7,0.9967567567567567)(8,0.9991891891891892)(9,0.9988513513513514)(10,0.9989864864864866)(11,0.999527027027027)(12,0.9983783783783783)(13,0.9978378378378379)(14,0.8580405405405405)(15,0.8592567567567566)(16,0.8572972972972971)(17,0.9991216216216215)(18,0.9996621621621622)};			
			
		\end{axis}
	\end{tikzpicture}
	
	\begin{tikzpicture}
		\begin{axis}[
			ybar,
			enlargelimits=0.12,
			legend style={at={(0.5,-0.55)}, anchor=north,legend columns=-1}, 
			ylabel={Node 2},
			symbolic x coords={1,2,3,4,5,6,7,8,9,10,11,12,13,14,15,16,17,18},
			xtick=data,
			nodes near coords align={vertical},
			bar width=5pt,
			]
			
			\addplot coordinates {(1,0.8976850040563794)(2,0.911762523275413)(3,0.9318227189688723)(4,0.9361804774332552)(5,0.847496529763766)(6,0.9283427692453724)(7,0.9244059866018546)(8,0.8544581709956649)(9,0.9787558248742546)(10,0.9194364478355862)(11,0.8630186588466813)(12,0.8747140507961436)(13,0.8521819498892126)(14,0.9538010687038511)(15,0.8147208940033476)(16,0.9233426651217517)(17,0.8516921145120737)(18,0.9158776206242458)};
			
			\addlegendentry{Weight}
			
			\addplot coordinates {(1,0.9934459459459457)(2,0.9958108108108115)(3,0.9983108108108105)(4,0.9939864864864862)(5,0.9933108108108109)(6,0.9960135135135134)(7,0.9899999999999998)(8,0.9985810810810808)(9,0.9949999999999997)(10,0.9979054054054056)(11,0.9968243243243247)(12,0.9931081081081081)(13,0.9903378378378378)(14,0.9924324324324322)(15,0.987837837837838)(16,0.9927027027027027)(17,0.9927027027027027)(18,0.9908108108108111)};		
			
			\addlegendentry{Volume}
		\end{axis}
		\draw [thin, <->] (rel axis cs:0.12,-0.2) -- (rel axis cs:1.11,-0.2) node[midway, below]
		{Front . . . . . . . . . . . . . . . . . . . . Pallets . . . . . . . . . . . . . . . . . . . . Rear};
		
	\end{tikzpicture}
	
	\caption{Occupation rates obtained by {\it Shims} with $scenario=1$, $surplus=1.2$\/ and $tmax=3600s$}
	\label{fig:packed2}
\end{figure}

Figure \ref{fig:packed2} shows the pallet occupancy rate (weight and volume) at each tour node found by {\it Shims}, with $scenario=1$, $surplus=1.2$\/ and $tmax =3600s$. As can be seen, the number of pallets with a high volume rate tends to grow throughout the tour.

For these reasons, we will present only the results obtained by {\it Shims}, with the average of the objective function $f$\/ and the runtime of {\it Gurobi}\/ and {\it Shims}. To facilitate the comparison between both, we added a last column in the tables where two values are indicated:

\begin{itemize}
	\item {\bf Normalized}: value between 0 and 1, which corresponds to the ratio between the sum of $f$\/ values obtained by the method in all scenarios and the sum of the best values obtained by both methods in all scenarios. The higher the value of {\bf Normalized}, the closer the method approached the best solutions found.
	\item {\bf Speed-up}: ratio of the sums of the runtimes of all scenarios and the sum of the method runtimes in all scenarios. The method with the highest {\bf Speed-up}\/ is the fastest.
\end{itemize}

We also indicate the adopted strategies: dedicating all the processing time to the $ntours=2$ shortest tours or distributing it among all $ntours=K!$\/ tours. The results obtained with $tmax = 3600s$, which is the highest tested runtime limit, are in Tables \ref{tab:20}, \ref{tab:50} and \ref{tab:100}, with $surplus$\/ values of $1.2$, $1.5$\/ and $2.0$, respectively. We indicate with an {\bf x} the cases where {\it Gurobi}\/ did not find a feasible solution within this runtime limit or had to be aborted due to high random-access memory (RAM) usage.

\begin{table}[htp] 
	\centering
	\caption{Solutions with $surplus = 1.2$\/ and $tmax =  3600s$}  \label{tab:20}
	\footnotesize
	\begin{tabular}{ccccccccc}
		\toprule
		$ntours$&$method$&$scenario$&{\bf 1}&{\bf 2}&{\bf 3}&{\bf 4}&{\bf 5}&\specialcell{{\bf Normalized}\\{\bf Speed-up}}   \\
		\toprule
		\multirow{4}{*}{2} &\multirow{2}{*}{{\it Gurobi}}  & $f$ & 8.53    & 11.79   & 13.14   & 13.52 & {\bf x}       &  0.9998 \\ 
		&                              &               {time (s)}   & 29      & 28      & 25      & 27    & {\bf x}          &  1.0    \\ 
		\cmidrule(lr){2-9}
		&                           \multirow{2}{*}{{\it Shims}}  & $f$ & 8.54    & 11.78   & 13.06   & 13.51 & 48.47       &  0.9980 \\
		&                             &                {time (s)}   & 1       & 1       & 1       & 1     & 2           &  22.7   \\
		
		\midrule
		
		\multirow{4}{*}{$K!$}&  \multirow{2}{*}{{\it Gurobi}} & $f$ & 8.60    & 12.20   & 13.66   & 15.00  & {\bf x}    & 0.9998 \\ 
		&                              &               {time (s)}   & 30      & 35      & 123     & 314    & {\bf x}    & 1.0 \\
		\cmidrule(lr){2-9}	
		&                           \multirow{2}{*}{{\it Shims}}  & $f$ & 8.49    & 12.10   & 13.30   & 14.49  & 53.61      & 0.9958  \\
		&                              &               {time (s)}   & 1       & 2       & 8       & 10     & 36         & 7.49   \\
		\bottomrule
		
	\end{tabular}
	\normalsize
\end{table}

\begin{table}[htp]
	\centering
	\caption{Solutions with $surplus = 1.5$\/ and $tmax =  3600s$}  \label{tab:50}
	\footnotesize
	\begin{tabular}{ccccccccc}
		\toprule
		$ntours$ & $method$          & $scenario$ & {\bf 1} & {\bf 2} & {\bf 3} & {\bf 4} & {\bf 5} & \specialcell{{\bf Normalized}\\{\bf Speed-up}}   \\
		\toprule
		\multirow{4}{*}{ 2 } & \multirow{2}{*}{{\it Gurobi}} & $f$ & 11.83 & 16.73 & 18.07 & 18.83 & 16.86 & 0.9996 \\ 
		&&               {time (s)}                                   & 55    & 64    & 39    & 40    & 88    & 1.0    \\ 
		\cmidrule(lr){2-9}	
		&\multirow{2}{*}{{\it Shims}}                               & $f$ & 11.85 & 16.72 & 18.05 & 18.80 & 16.87 & 0.9993  \\
		&&         {time (s)}                                         & 1     & 1     & 2     & 2     & 2     & 35.8  \\
		
		\midrule
		
		\multirow{4}{*}{$K!$} & \multirow{2}{*}{{\it Gurobi}}   & $f$ & 11.83 & 16.93 & 18.40 & 20.95 & 17.60 & 0.9999 \\  
		&&               {time (s)}                                   & 63    & 59    & 195   & 472   & 2,258  & 1.0  \\
		\cmidrule(lr){2-9}	
		&\multirow{2}{*}{{\it Shims}}                               & $f$ & 11.85 & 16.91 & 18.36 & 20.93 & 17.50 & 0.9976  \\ 
		&&         {time (s)}                                         & 1     & 2     & 5     & 15    & 100   & 23.8 \\
		\bottomrule
		
	\end{tabular}
	\normalsize
\end{table}

\begin{table}[htp]
	\centering
	\caption{Solutions with $surplus = 2.0$\/ and $tmax =  3600s$}  \label{tab:100}
	\footnotesize
	\begin{tabular}{ccccccccc}
		\toprule
		$ntours$& $method$& $scenario$ & {\bf 1} & {\bf 2} & {\bf 3} & {\bf 4} & {\bf 5} & \specialcell{{\bf Normalized}\\{\bf Speed-up}}   \\
		\toprule
		\multirow{4}{*}{2} & \multirow{2}{*}{{\it Gurobi}} & $f$ & 17.70 & 24.20 & 26.39 & 27.17 & 24.20 & 0.9995 \\ 
		&&                                              {time (s)}  & 168   & 98    & 79    & 70    & 72    & 1.0  \\
		\cmidrule(lr){2-9}	
		&                            \multirow{2}{*}{{\it Shims}} & $f$ & 17.74 & 24.22 & 26.32 & 27.07 & 23.13 & 0.9896  \\
		&&                                              {time (s)}  & 1     & 2     & 2     & 3     & 4     & 40.6  \\
		
		\midrule
		
		\multirow{4}{*}{$K!$} & \multirow{2}{*}{{\it Gurobi}} & $f$ & 17.90 & 25.44 & 26.51 & 29.13 &{\bf x}& 0.9994 \\ 
		&&                                               {time (s)} & 178   & 143   & 378   & 862   &{\bf x}& 1.0  \\
		\cmidrule(lr){2-9}	
		&                         \multirow{2}{*}{{\it Shims}}    & $f$ & 17.94 & 25.45 & 26.44 & 28.84 & 26.22  & 0.9970  \\
		&&                                               {time (s)} & 1     & 3     & 10    & 31    & 196    & 34.7  \\
		\bottomrule
		
	\end{tabular}
	\normalsize
\end{table}

From these data, we can draw some conclusions:

\begin{itemize}

	\item The strategy of testing all $K!$\/ tours often provide a better-quality solution, even with less time on each node. This shows that the four sub-problems are interconnected in such a way that it is not enough to solve them separately.
	
	\item {\it Gurobi}\/ fails in some cases when $scenario=5$\/ and the strategy is to check all $K!$\/ tours. This occurs because the runtime limit per node is smaller and there tend to be more packed contents on the aircraft, reducing the space for allocating items and making the solution difficult.
	
	\item When {\it Gurobi}\/ finishes, it finds the best solution, but the one obtained by {\it Shims}\/ reaches at least $98.96\%$\/ of that value. Considering only the strategy of testing all $K!$ tours, this value increases to $99.58\%$.
	
	\item {\it Shims}\/ always finds a solution, being 7 to 40 times faster.
	
	\item All runtimes are much lower than the limit because the solution on many nodes can be fast. Anyway, in all the tests performed, the maximum time spent by {\it Shims}\/ did not reach 4 minutes. On the other hand, when $scenario=5$\/ and $surplus=1.5$, {\it Gurobi}\/ spent almost 40 minutes.
	
\end{itemize}

Table \ref{tab:100time} shows the results obtained with the strategy of testing the $K!$\/ tours in all scenarios with different $tmax$. We can observe more cases where {\it Gurobi}\/ fails, even in smaller scenarios. When {\it Gurobi}\/ finishes, {\it Shims}\/ finds a solution of similar quality ($99\%$\/ or better). In all cases, {\it Shims}\/ finds a solution in less than 4 minutes.

\begin{table}[H]
	\centering
	\caption{Solutions testing all $K!$\/ tours with different runtime limits}  \label{tab:100time}
	\scriptsize
	\setlength{\tabcolsep}{3.5pt}
	\begin{tabular}{ccc|ccccc|ccccc|ccccc}
		\toprule
		&             \multicolumn{2}{r}{$surplus$} &\multicolumn{5}{c}{\bf 1.2}&\multicolumn{5}{c}{\bf 1.5}&\multicolumn{5}{c}{\bf 2.0} \\
		\midrule
		$method$& {$tmax$} & $scenario$  &{\bf 1}&{\bf 2}&{\bf 3}&{\bf 4}&{\bf 5}&{\bf 1}&{\bf 2}&{\bf 3}&{\bf 4}&{\bf 5}&{\bf 1}&{\bf 2}&{\bf 3}&{\bf 4}&{\bf 5} \\
		\toprule
		
		\multirow{10}{*}{{\it Gurobi}} & \multirow{2}{*}{240s} & $f$ & 8.60 & 12.20 & 13.67 & {\bf x} & {\bf x} & 11.77 & 16.38 & 18.10 & {\bf x} & {\bf x} & 17.89 & 25.42 & {\bf x} & {\bf x} & {\bf x} \\ 
		&               &  {time (s)}                            & 31   & 35    & 124   & {\bf x} & {\bf x} & 52    & 59    & 200   & {\bf x} & {\bf x} & 188   & 145   & {\bf x} & {\bf x} & {\bf x} \\ 
		\cmidrule(lr){2-18}		
		
		&                           \multirow{2}{*}{1200s} & $f$ & 8.61 & 12.20 & 13.31 & 10.00   & {\bf x} & 11.77 & 17.02 & 18.25 &  20.64  & {\bf x} & 17.75 & 25.24 & 26.49 & 27.97   & {\bf x} \\ 
		&        &  {time (s)}                                   & 28   & 37    & 129   & 320     & {\bf x} & 46    & 61    & 190   &  304    & {\bf x} & 161   & 139   & 384   & 579     & {\bf x} \\ 
		\cmidrule(lr){2-18}                              
		&                           \multirow{2}{*}{2400s} & $f$ & 8.60 & 12.21 & 13.67 & 15.00   & 13.41   & 11.77 & 16.37 & 18.03 & 20.95   & {\bf x} & 17.89 & 25.44 & 26.15 & 29.13   & {\bf x} \\
		&        &  {time (s)}                                   & 26   & 38    & 134   & 310     & 1,520    & 46    & 60    & 199   & 461     & {\bf x} & 164   & 140   & 383   & 786     & {\bf x} \\
		\cmidrule(lr){2-18}                              		
		
		&                           \multirow{2}{*}{3600s} & $f$ & 8.60 & 12.20 & 13.66 & 15.00   & {\bf x} & 11.76 & 16.37 & 18.01 &  20.95  & 17.60   & 17.90 & 25.44 & 26.15 & 29.13   & {\bf x} \\ 
		&        &  {time (s)}                                   & 30   & 35    & 123   & 314     & {\bf x} & 64    & 58    & 195   &  472    & 2,258    & 178   & 143   & 378   & 862     & {\bf x} \\ 
		
		\midrule
		
		\multirow{2}{*}{{\it Shims}} &  \multirow{2}{*}{240s}  & $f$ & 8.49 & 12.10 & 13.30 & 14.49   & 53.61   & 11.78 & 16.29 & 17.99 &  20.93  & 17.50   & 17.94 & 25.45 & 26.14 & 28.84   & 26.22 \\ %
		&         & {time (s)}                                       & 1    & 2     & 8     & 10      & 36      & 1     & 2     & 5     &  15     & 100     & 1     & 3     & 10    & 31      & 196   \\ 
		
		\bottomrule
	\end{tabular}
	\normalsize
\end{table}

The actual RAM consumption of {\it Gurobi}\/ was over 8.5 GB, and all of {\it Shims}\/'s executions consumed at most 1.5 GB of RAM.

\subsubsection{Results when $K>6$}

These last results do not correspond to practical cases of air transport, as tours where $K>6$ very rarely occur. However, it is possible to see that the {\it Shims}\/ maintains robust behaviour as the number of nodes grows, that is, it could be adapted to similar contexts (ships and trucks, for example), where there may be more nodes.

Considering real data from the 15 main Brazilian airports, we implemented a GA-based TSP heuristic that returned 100 tours in approximately $33 s$. We implemented this heuristic with DEAP (Distributed Evolutionary Algorithms in Python), an evolutionary computation framework. For more details, see \citep{DEAP_JMLR2012} and {\tt github.com/deap/deap}. 

Figure \ref{fig:performance} shows the runtime curve of {\it Shims}\/ as the number $K$\/ of nodes increases. Runtime is the average obtained from 7 instances generated with $surplus=2.0$\/ and $tmax=1200s$\/ for each value of $K$. In Figure \ref{fig:morenodes}, we indicate one of the tours found by this TSP heuristic when $K=15$.

\begin{table}[H]
	
\begin{minipage}{0.58\linewidth}
		\centering
		\footnotesize
		\pgfplotsset{compat=newest}

\begin{tikzpicture}
	\begin{axis}[
		xlabel style={align=center},
		xlabel={Number $K$\/ of nodes},
		ylabel={Runtime (s)},
		xmin=7,   
		ymin=200, 
		xtick={ 5, 6, 7,  8,  9, 10, 11, 12, 13, 14, 15, 16},
		ytick={250,350,450,550,650,750,850,950,1050,1150,1250},	
		ymajorgrids=true,
		grid style=dashed,
		]
		\addplot[
		color=blue,
		mark=square,
		]
		coordinates {
			(7,273)(8,360)(9,468)(10,560)(11,660)(12,756)(13,884)(14,1008)(15,1170)
		};
		
	\end{axis}
\end{tikzpicture}
	
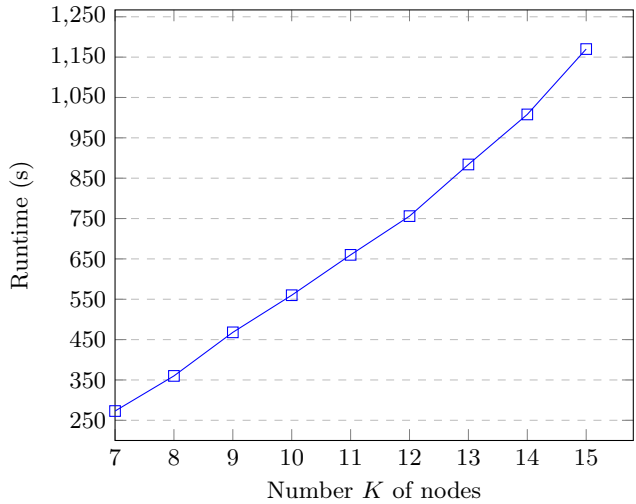
\captionof{figure}{{\it Shims}\/ performance with $surplus=2.0$ and $tmax = 1200s$}
\label{fig:performance}
\end{minipage}\hfill 
\begin{minipage}{0.42\linewidth}
\centering
\includegraphics[scale=0.293]{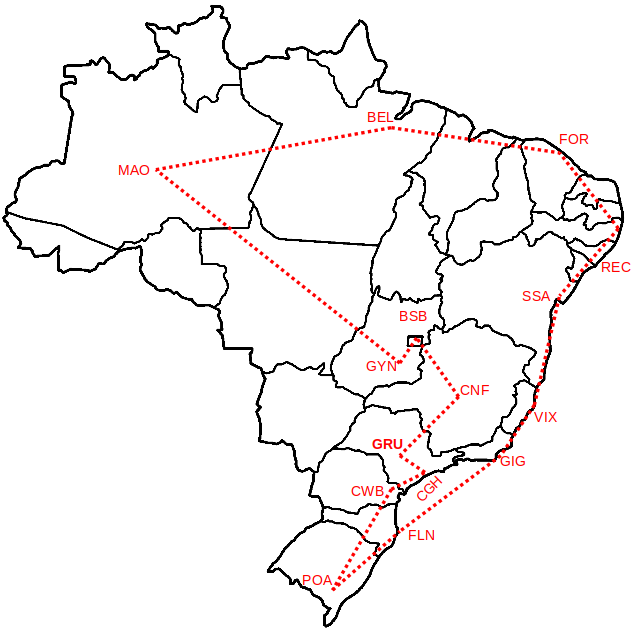}
\captionof{figure}{A tour where $K=15$}
\label{fig:morenodes}
\end{minipage}
\end{table}

\section{Conclusions}
\label{conclusions}

In this work, we modelled and solved a real air transport problem named {\it Air Cargo Load Planning with Routing, Pickup, and Delivery Problem}\/ (ACLP+RPDP). For the first time in the literature, a {\it NP-hard}\/ problem that involves {\it simultaneously}\/ pallet assembly, load balancing, route planning, and pickup and delivery is addressed, where the cost-effectiveness of transport is maximized. Currently, there is no commercial software available for this problem.

We adopted some simplifications that are not critical, but that allowed for an unprecedented solution to this problem considering several nodes. In practical cases, there are hundreds of items to be shipped at each node, and the number of nodes, excluding the base, is smaller than the number of pallets. Considering a real aircraft, we have developed node-by-node solutions such that the complete process can be executed quickly on a handheld computer, offering good results and reducing stress for the transport planners.

As validation, we carried out tests in several scenarios. In the real cases, the solution process can establish, in less than four minutes, a flight itinerary for a single aircraft with a good distribution of load on pallets at each node of the tour, enforcing the weight balance, maximising the total score, and minimising fuel consumption along the planned route, which is beneficial to reduce carbon emissions. This output is an essential part of airlift: it guarantees flight safety, makes ground operations more efficient, and makes sure that each item gets to its right destination.

Our main contributions were the mathematical modelling of ACLP+RPDP, involving four well-known and interconnected {\it NP-hard}\/ sub-problems, a complete process to solve it, and a new heuristic named {\it Shims}\/ that offers fast node solutions with good quality. Without a quick method for balanced allocation at each node, it would be unfeasible to find a flight itinerary and, consequently, a pickup and delivery plan on this tour.

We also show that this solution process remains valid in other contexts where there may be a greater number of nodes. This way, our method is not exclusive to aircraft and airports: it can be adapted to ships and ports, trucks and warehouses, or wagons and railways. In these situations, it would be necessary to make some changes in the model: for example, modify the load balancing constraints and consider parts of the available space as pallets.

As this is ongoing research, we thought about some possible future improvements: consider more than one aircraft, implement parallel algorithms in some steps of the solution to improve computational efficiency, and model 3-D items.

\section*{CRediT authorship contribution statement}

\textbf{Antonio Celio Pereira de Mesquita:} Conceptualization, Methodology, Software, Writing - original draft preparation, Investigation, Validation. \textbf{Carlos Alberto Alonso Sanches:} Conceptualization, Methodology, Resources, Supervision, Writing - reviewing \& editing.

\section*{Declaration of competing interest}

The authors declare that they have no known competing financial interests or personal relationships that could have appeared to influence the work reported in this paper.

\section*{Data availability}

The dataset and code used in this work are available in {\tt github.com/celiomesquita/ACLP\_RPDP\_P}.

\section*{Acknowledgments}

This research was partially supported by the S\~{a}o Paulo Research Foundation (FAPESP), grant 2022/05803-3.

\section*{References}

\bibliographystyle{apalike} 

\bibliography{bib}

@ARTICLE{Bertsimas,
	author = {Bertsimas, Dimitris and Chang, Allison and Misic, Velibor V. and Mundru, Nishanth},
	title = {The airlift planning problem},
	year = {2019},
	journal = {Transportation Science},
	volume = {53},
	number = {3},
	pages = {773-795},
	doi = {10.1287/trsc.2018.0847},
	type = {Article}
}

@InProceedings{Golestanian,
	author =	{Golestanian, Arnoosh and Bianco, Giovanni Lo and Tao, Chengyu and Beck, J. Christopher},
	title =	{{Optimization Models for Pickup-And-Delivery Problems with Reconfigurable Capacities}},
	booktitle =	{29th International Conference on Principles and Practice of Constraint Programming (CP 2023)},
	pages =	{17:1--17:17},
	series =	{Leibniz International Proceedings in Informatics (LIPIcs)},
	ISBN =	{978-3-95977-300-3},
	ISSN =	{1868-8969},
	year =	{2023},
	volume =	{280},
	editor =	{Yap, Roland H. C.},
	publisher =	{Schloss Dagstuhl -- Leibniz-Zentrum f{\"u}r Informatik},
	address =	{Dagstuhl, Germany}
}

@ARTICLE{Meng,
	author = {Meng, Shanshan and Guo, Xiuping and Li, Dong and Liu, Guoquan},
	title = {The multi-visit drone routing problem for pickup and delivery services},
	year = {2023},
	journal = {Transportation Research Part E: Logistics and Transportation Review},
	volume = {169},
	pages = {102990},
	doi = {10.1016/j.tre.2022.102990}
}

@ARTICLE{Debnath,
	author = {Debnath, Dipraj and Hawary, A.F.},
	title = {Adapting Travelling Salesmen Problem for Real-Time \mbox{UAS} Path Planning Using Genetic Algorithm},
	year = {2021},
	journal = {Lecture Notes in Mechanical Engineering},
	pages = {151-163},
	doi = {10.1007/978-981-16-0866-7_12},
	type = {Conference paper}
}

@article{Cheikh,
	title={A comprehensive survey on the Multiple Traveling Salesman Problem: Applications, approaches and taxonomy},
	volume={40},
	ISSN={1574-0137},
	url={http://dx.doi.org/10.1016/j.cosrev.2021.100369},
	DOI={10.1016/j.cosrev.2021.100369},
	journal={Computer Science Review},
	publisher={Elsevier BV},
	author={Cheikhrouhou, Omar and Khoufi, Ines},
	year={2021},
	month=may, pages={100369} }

@ARTICLE{Xie,
	author={Xie, Junfei and Carrillo, Luis Rodolfo Garcia and Jin, Lei},
	journal={IEEE Control Systems Letters}, 
	title={An Integrated Traveling Salesman and Coverage Path Planning Problem for Unmanned Aircraft Systems}, 
	year={2019},
	volume={3},
	number={1},
	pages={67-72},
	doi={10.1109/LCSYS.2018.2851661}}

@INPROCEEDINGS{Ahmad,
	author={Ahmad, Edwin Dwi and Muklason, Ahmad and Nurkasanah, Ika},
	booktitle={2020 International Conference on Computer Engineering, Network, and Intelligent Multimedia (CENIM)}, 
	title={Route Optimization of Airplane Travel Plans Using the Tabu-Simulated Annealing Algorithm to Solve the Traveling Salesman Challenge 2.0}, 
	year={2020},
	volume={},
	number={},
	pages={217-221},
	doi={10.1109/CENIM51130.2020.9297892}}

@article{DEAP_JMLR2012, 
	author    = " F\'elix-Antoine Fortin and Fran\c{c}ois-Michel {De Rainville} and Marc-Andr\'e Gardner and Marc Parizeau and Christian Gagn\'e ",
	title     = { {DEAP}: Evolutionary Algorithms Made Easy },
	pages    = { 2171--2175 },
	volume    = { 13 },
	month     = { jul },
	year      = { 2012 },
	journal   = { Journal of Machine Learning Research }
}

@inproceedings{PeerlinckAmy2022MFEO,
	pages = {1--8},
	publisher = {IEEE},
	booktitle = {2022 IEEE Congress on Evolutionary Computation (CEC)},
	isbn = {9781665467087},
	year = {2022},
	title = {Multi-Objective Factored Evolutionary Optimization and the Multi-Objective Knapsack Problem},
	language = {eng},
	author = {Peerlinck, Amy and Sheppard, John},
}

@article{MiguelOptimalAP,
	title={Optimal Aircraft Payload Weight and Balance using Fuzzy Linear Programming Model},
	author={Juan Miguel V. Macalintal and Aristotle T. Ubando},
	volume={103},
	year={2023},	
	journal={Chemical Engineering Transactions},
	pages={613-618},
	doi={DOI:10.3303/CET23103103},
}

@article{LopezIbanezManuel2016,
	journal = {Operations {R}esearch {P}erspectives},
	pages = {43--58},
	volume = {3},
	year = {2016},
	title = {The irace package: iterated racing for automatic algorithm configuration},
	language = {eng},
	author = {M. Lopez-Ibanez and J. Dubois-Lacoste and L.P. C\'{a}ceres and M. Birattari and T. St\"{u}tzle},	
}

@article{JohnsonGarey1985,
	issn = {0885-064X},
	journal = {Journal of Complexity},
	pages = {65--106},
	volume = {1},
	number = {1},
	year = {1985},
	title = {A 7160 theorem for bin packing},
	language = {eng},
	author = {Johnson, David S and Garey, Michael R},
}

@inproceedings{fok2004optimizing,
	title={Optimizing air cargo load planning and analysis},
	author={Fok, Kelly and Chun, A},
	booktitle={Proceedings of the International Conference on Computing, Communications and Control Technologies},
	year={2004},
	pages = {520--531},
}

@article{Kevin1992,
	author = {Kevin Y. K. Ng},
	year = {1992},
	title = {A Multicriteria Optimization Approach to Aircraft Loading},
	journal = {Operations Research},
	volume = {40},
	number = {6},
	pages = {1200-1205},
	publisher = {INFORMS},
}

@article{LarsenMikkelsen1979,
	author = {Ole Larsen and Gert Mikkelsen},
	title = {An interactive system for the loading of cargo aircraft},
	journal = {European Journal of Operational Research},
	year = {1980},
	volume = {4},
	pages = {367-373},
	number = {6},
	publisher = {Elsevier},
}

@article{eugene2021,
	author = {Eugene Y. C. Wong and Daniel Y. Mo and Stuart So},
	title = {Closed-loop digital twin system for air cargo load planning operations},
	journal = {International Journal of Computer Integrated Manufacturing},
	volume = {34},
	number = {7-8},
	pages = {801-813},
	year  = {2021},
	publisher = {Taylor and Francis},
	doi = {10.1080/0951192X.2020.1775299},
}

@INPROCEEDINGS{wong2020,
	author={Wong, Eugene Y.C. and Ling, Kev K. T.},
	booktitle={7th International Conference on Frontiers of Industrial Engineering (ICFIE)}, 
	title={A Mixed Integer Programming Approach to Air Cargo Load Planning with Multiple Aircraft Configurations and Dangerous Goods}, 
	year={2020},
	pages={123-130},
	doi={10.1109/ICFIE50845.2020.9266727},
}

@article{zhao2021,
	author = {Zhao, Xiangling and Yuan, Yuan and Dong, Yun and Zhao, Ren},
	title = {Optimization approach to the aircraft weight and balance problem with the centre of gravity envelope constraints},
	journal = {IET Intelligent Transport Systems},
	volume = {15},
	number = {10},
	pages = {1269-1286},
	publisher = {IET},
	year = {2021},
}

@article{zhao2023,
	author = {X. Zhao and Y. Dong and L. Zuo},
	title = {A combinatorial optimization approach for air cargo palletization and aircraft loading},
	journal = {Mathematics},
	volume = {11},
	number = {13},
	pages = {1-16},
	year = {2023},
}

@article{BrandtStefan2019,
	issn = {0377-2217},
	journal = {{European Journal of Operational Research}},
	pages = {399-410},
	volume = {275},
	publisher = {Elsevier},
	number = {2},
	year = {2019},
	title = {The air cargo load planning problem - a consolidated problem definition and literature review on related problems},
	author = {Brandt, Felix and Nickel, Stefan},
}

@article{Brosh1981,
	author = {Israel Brosh},
	title = {Optimal cargo allocation on board a plane: a sequential linear programming approach},
	journal = {European Journal of Operational Research},
	year = {1981},
	pages = {40-46},
	volume = {8},
	number = {1},
	publisher = {Elsevier},
}

@article{Chan2006,
	author = {F.T.S. Chan and R. Bhagwat and N. Kumar and M.K. Tiwari and P. Lam},
	title = {Development of a decision support system for air-cargo pallets loading problem: A case study},
	journal = {Expert Systems with Applications},
	year = {2006},
	volume = {31},
	number = {3},
	pages = {472-485},
	publisher = {Elsevier},
}

@article{CharonHudry2001,
	issn = {0377-2217},
	journal = {European Journal of Operational Research},
	pages = {86-101},
	volume = {135},
	publisher = {Elsevier},
	number = {1},
	year = {2001},
	title = {The noising methods: A generalization of some metaheuristics},
	author = {Charon, Irène and Hudry, Olivier},
}

@article{CharonHudry1993,
	issn = {0167-6377},
	journal = {Operations Research Letters},
	pages = {133-137},
	volume = {14},
	publisher = {Elsevier},
	number = {3},
	year = {1993},
	title = {The noising method: a new method for combinatorial optimization},
	author = {Charon, Irène and Hudry, Olivier},
}

@article{DorigoManiezzoColorni1996,
	author = {M. Dorigo and V. Maniezzo and A. Colorni},
	title = {The ant system: optimization by a colony of cooperating agents},
	journal = {IEEE Transactions on Systems, Man, and Cybernetics},
	year = {1996},
	volume = {26},
	pages = {29-41},
}

@phdthesis{Dorigo1992,
	author  = {M. Dorigo},
	title   = {Optimization, Learning and Natural Algorithms},
	school  = {Politecnico di Milano},
	year    = {1992}
}

@article{FeoResende1989,
	author = {T. A. Feo and M. G. C. Resende},
	title = {A probabilistic heuristic for a computationally difficult set covering problem},
	journal = {Operations Research Letters},
	year = {1989},
	volume = {8},
	pages = {67-71},
}

@article{Glover1986,
	author = {F. Glover},
	title = {Future paths for integer programming and links to artificial intelligence},
	journal = {Computers and Operations Research},
	year = {1986},
	volume = {13},
	pages = {533-549},
}

@article{Heidelberg1998,
	author = {Kurt R. Heidelberg and Gregory S. Parnell and James E. Ames},
	title = {Automated air load planning},
	journal = {Naval Research Logistics},
	year = {1998},
	volume = {45},
	number = {8},
	pages = {751-768},
	publisher = {John Wiley \& Sons},
}

@BOOK{Fidanova2006,
	booktitle = {Handbook of Research on Nature Inspired Computining for Economics and Management},
	publisher = {J--Ph. Renard editor},
	year = {2006},
	author = {Stefka Fidanova},
	title = {Ant Colony Optimization and Multiple Knapsack Problem},
	volume = {Chapter 33},
	institution={Idea Group Inc},
	key = {ISBN 1--59140-984--5},
	pages={498-509},
}

@article{Limbourg2012,
	author = {S. Limbourg and M. Schyns and G. Laporte},
	title = {Automatic aircraft cargo load planning},
	journal = {Journal of the Operational Research Society},
	year = {2012},
	volume = {63},
	number = {9},
	pages = {1271-1283},
}

@article{LurkinSchyns2015,
	author = {Lurkin, Virginie and Schyns, Michaël},
	title = {The Airline Container Loading Problem with pickup and delivery},
	journal = {{European Journal of Operational Research}},
	year = {2015},
	volume = {244(3)},
	pages = {955-965},
	publisher = {Elsevier},
}

@article{MongeauBes2003,
	author = {M. Mongeau and C. Bes},
	title = {Optimization of aircraft container loading},
	journal = {IEEE Transaction on Aerospace and Electronic Systems},
	year = {2003},
	volume = {39},
	number = {1},
	pages = {140-150},
	publisher = {IEEE},
}

@INPROCEEDINGS{NiarFreville1997,
	author = {S. Niar and A. Freville},
	title = {A parallel tabu search algorithm for the 0-1 multidimensional knapsack problem},
	booktitle = {Proceedings 11th International Parallel Processing Symposium},
	year = {1997},
	pages = {512-516},
}

@article{Paquay2016,
	author = {Paquay, C. and Schyns, M. and Limbourg, S.},
	title = {A mixed integer programming formulation for the three-dimensional bin packing problem deriving from an air cargo application},
	journal = {International Transactions in Operational Research},
	year = {2016},
	volume = {23},
	pages = {187-213},
}

@article{RoesenerHall2014,
	title={A NONLINEAR INTEGER PROGRAMMING FORMULATION FOR THE AIRLIFT LOADING PROBLEM WITH INSUFFICIENT AIRCRAFT},
	volume={5},
	number={1},
	journal={Journal of Nonlinear Analysis and Optimization: Theory and Applications},
	author={Roesener, A. and Hall, S.},
	year={2014},
	pages={125-141} 
}

@article{RoesenerBarnes2016,
	author = {A.G. Roesener and J.W. Barnes},
	title = {An advanced tabu search approach to the dynamic airlift loading problem},
	journal = {Logistics Research},
	year = {2016},
	volume = {9(1)},
	pages = {1-18},
}

@article{Vancroonemburg2014,
	issn = {1366-5545},
	journal = {Transportation Research Part E: Logistics and Transportation Review},
	pages = {70-83},
	volume = {65},
	publisher = {Elsevier},
	year = {2014},
	title = {Automatic air cargo selection and weight balancing: A mixed integer programming approach},
	author = {Vancroonenburg, Wim and Verstichel, Jannes and Tavernier, Karel and Vanden Berghe, Greet},
}

@article{Verstichel2011,
	author = {J. Verstichel and W. Vancroonenburg and W. Souffriau and G. V. Berghe},
	title = {A mixed integer programming approach to the aircraft weight and balance problem},
	journal = {Procedia Social and Behavioral Sciences},
	year = {2011},
	volume = {20},
	pages = {1051-1059},
	publisher = {Elsevier},
}

@article{Zhan2020,
	author = {S. Zhan and L. Wang and Z. Zhang and Y. Zhong},
	title = {Noising methods with hybrid greedy repair operator for 0-1 knapsack problem},
	journal = {Memetic Computing},
	year = {2020},
	volume = {12},
	pages = {37-50},
}

@article{Alonso2019,
	author = {M. T. Alonso and R. Alvarez-Valdes and F. Parreno},
	title = {A \mbox{GRASP} algorithm for multi-container loading problems with practical constraints},
	journal = {A Quarterly Journal of Operations Research},
	year = {2019},
	volume = {18},
	pages = {49-72},
}

@article{KaluznyBohdanL2009Oalb,
	issn = {0969-6016},
	journal = {International Transactions in Operational Research},
	pages = {767-787},
	volume = {16},
	publisher = {John Wiley \& Sons},
	number = {6},
	year = {2009},
	title = {Optimal aircraft load balancing},
	author = {Kaluzny, Bohdan L and Shaw, R. H. A. David},
}

@article{YangLiuGao2018,
	issn = {2261-236X},
	journal = {MATEC Web of Conferences},
	pages = {1-6},
	volume = {179},
	publisher = {EDP Sciences},
	year = {2018},
	title = {Load Planning of Transport Aircraft Based on Hybrid Genetic Algorithm},
	author = {Chenguang, Yang and Hu, Liu and Yuan, Gao},
}

@article{Paquay2018,
	author = {Paquay, C. and Schyns, M. and Limbourg, S. and Oliveira, J. F.},
	title = {{MIP}-based constructive heuristics for the three-dimensional {Bin Packing Problem} with transportation constraints},
	journal = {International Journal of Production Research},
	volume = {56},
	number = {4},
	pages = {1581-1592},
	year  = {2018},
	publisher = {Taylor and Francis},
}

@book{holland1992adaptation,
	title={Adaptation in natural and artificial systems: an introductory analysis with applications to biology, control, and artificial intelligence},
	author={Holland, John H},
	year={1992},
	publisher={MIT press}
}

\end{document}